\documentclass[10pt,reqno]{article}
\usepackage{fullpage}

\usepackage[unicode=true]{hyperref}
\usepackage{amsmath}
\usepackage{cite}
\usepackage{amsfonts}
\usepackage{amssymb}
\usepackage{amsthm}
\usepackage{authblk}

\usepackage[pdftex]{color,graphicx}
\usepackage{mypack} 
\usepackage{adjustbox}
\usepackage{soul}
\usepackage{comment}
\usepackage{bbold}
\usepackage{tikz}
\usetikzlibrary{decorations.pathreplacing,angles,quotes}
\usepackage{bbold}

\usepackage{diagbox}
\usepackage{enumerate}

\usetikzlibrary{matrix,arrows,decorations.pathmorphing}
\usepackage{tikz-cd}
\usetikzlibrary{arrows} 

\usetikzlibrary{decorations.markings}
 \usepackage{caption}
\usepackage{subcaption}
 
 \usepackage{pdfpages}

\usepackage{blkarray}

\usepackage{centernot}
\usepackage{mathtools}
\usepackage{stmaryrd}

  \usepackage{arydshln}

\newcommand{\zd}{\ensuremath{\mathbb{Z}_d}}

\definecolor{bluegray}{rgb}{0.4, 0.6, 0.8}
\definecolor{turquoise}{rgb}{0.2, 0.7, 0.6}
\definecolor{hy-green}{rgb}{0.1, 0.5, 0.1}
 
\usepackage{soul}  
 \newcommand{\suchthat}{\;\ifnum\currentgrouptype=16 \middle\fi|\;}

\title{Simplicial techniques for operator solutions of linear constraint systems}

\author[1]{Ho Yiu Chung}
\author[2]{Cihan Okay}
\author[2]{Igor Sikora}
\affil[1]{{\small{Department of Computer Science, The University of Hong Kong, Pokfulam, Hong Kong}}}
\affil[2]{{\small{Department of Mathematics, Bilkent University, Ankara, Turkey}}}

\begin{document}
  \maketitle

\begin{abstract}
{A linear constraint system is specified by linear equations over the group $\ZZ_d$ of integers modulo $d$. Their operator solutions play an important role in the study of quantum contextuality and non-local games.}
In this paper, we use the theory of simplicial sets to develop a framework for studying operator solutions of linear systems.
Our approach refines the well-known group-theoretical approach based on solution groups by identifying these groups as algebraic invariants closely related to the fundamental group of a space.
In this respect, our approach also makes a connection to the earlier homotopical approach based on cell complexes. Within our framework, we introduce a new class of linear systems that come from simplicial sets and show that any linear system can be reduced to one of that form.
Then we specialize in linear systems that are associated with groups. We provide significant evidence for a conjecture stating that for odd $d$ every linear system admitting a solution in a group admits a solution in $\ZZ_d$.
\end{abstract}

\tableofcontents

\section{Introduction}  

Linear (constraint) systems are a source of contextual distributions that arise in quantum theory \cite{Mer93,peres1991two,arkhipov2012extending} and play a prominent role in studying non-local games \cite{cleve2017perfect,slofstra2020tsirelson}.
A linear system over the group $\ZZ_d$ of integers modulo $d$ is specified by an equation $Ax=b$ where $A$ is an $r\times c$-matrix   and $b$ is a column of size $r$ both with entries in $\ZZ_d$.
A solution of such a linear system in a group $G$  
consists of group elements $T_1,\cdots,T_c$ satisfying 
$$
T_1^{A_{i1}} T_2^{A_{i2}} \cdots T_c^{A_{ic}} = J_G^{b_i}
$$ 
where $J_G$ is a fixed central element of order $d$ in the group.
In addition to these product equations, a solution has to satisfy (1) $d$-torsion property, that is $T_i^d$ for all $i=1,\cdots,c$, and (2) commutativity property: $T_iT_j=T_jT_i$ whenever $A_{ki}$ and $A_{kj}$ both non-zero for some row index $k$.
In the literature, it is common practice to take $G$ as a unitary group acting on a finite-dimensional Hilbert space. In this case, the solutions are usually called operator solutions.
There is a group-theoretic approach for
studying operator solutions of linear systems
centered around the properties of the solution group $\Gamma(A,b)$; see, for example, \cite{cleve2017perfect}. 
In this paper, we associate spaces to linear systems such that their algebraic invariants are closely related to the solution group.
Our methods connect both to the solution group via a construction akin to the fundamental group of a space, and the homotopical methods of \cite{okay2020homotopical} based on cell complexes.

{We use simpicial sets \cite{GJ99} as combinatorial models of spaces.}
Simplicial sets are fundamental objects of modern homotopy theory. They are more expressive than their close relatives simplicial complexes.
One can associate a simplicial complex $\Sigma$ to a matrix $A$ specifying a linear system:
The vertices of $\Sigma$ are given by $v_1,\cdots,v_c$ and maximal simplices by $\sigma_1,\cdots,\sigma_r$ where each $\sigma_i$ consists of vertices $v_j$ such that $A_{ij}\neq 0$. 
Each row of $A$ can be regarded as a function $A_i:\Sigma_0\to \ZZ_d$ on the {vertex set.}
A simplicial set consists of a set of $n$-simplices for each $n\geq 0$ together with the simplicial relations  describing how simplices of various dimensions are glued.  
Our simplicial realizations are motivated by a well-known construction in algebraic topology.
The nerve space $NG$ of a group consists of $n$-tuples of group elements as its $n$-simplices.
Given a matrix $A$ with the associated simplicial complex $\Sigma$ we construct a simplicial set, denoted by $N(\ZZ_d,\Sigma)$, whose set of $n$-simplices are given by tuples
$
(s_1,s_2,\cdots,s_n)
$ 
of functions $s_i:\Sigma_0\to \ZZ_d$ such that the union of the supports $\supp(s_i)$ is a simplex of $\Sigma$.
The rows $A_i$ can be regarded as $1$-simplices of this simplicial set.
Another construction we need is the  simplicial set $N(\ZZ_d,G)$, which we refer to as the $d$-torsion commutative nerve of $G$,
consisting of $n$-simplices given by tuples of pairwise commuting $d$-torsion group elements.
Closely related versions of the nerve construction are first introduced in \cite{ACT12}, and their homotopy theory has been studied recently; see, for example, \cite{O16,antolin2020classifying,okay2021commutative}.
 
\begin{custompro}{\ref{pro:solution}}
There is a bijective correspondence between solutions of $(A,b)$ in a group $G$ and the simplicial set maps $f$ making the following diagram commute
$$
\begin{tikzcd}[column sep=huge,row sep=large]
\vee_{i=1}^r N\ZZ_d \arrow[r,"\beta"] \arrow[d,"\alpha"',hook]  & N\ZZ_d \arrow[d,"\iota",hook] \\
N(\ZZ_d,\Sigma) \arrow[r,"f"] & N(\ZZ_d,G)
\end{tikzcd}
$$
\end{custompro}  

The maps that appear in this diagram are described in Section \ref{sec:SimpReal}. Briefly, $\alpha$, $\beta$ and $\iota$ correspond to $A_i$'s, $b_i$'s, and the central element $J_G$ in a way that the commutativity of the diagram coincides with the notion of a solution in $G$ introduced above.
The cofiber of the $\alpha$ map, which we denote by $\bar N(\ZZ_d,\Sigma)$, together with a cohomology class $\gamma_b$ serves as the simplicial realization of the linear system $(A,b)$. The cohomology class depends on the column vector $b$ and is defined using the cohomology long exact sequence associated to the cofiber sequence.
There is a converse to this procedure which associates a linear system $(A_X,b_\gamma)$ to a pair $(X,\gamma)$ {consisting} of a simplicial set $X$ and a cohomology class $[\gamma]\in H^2(X)$, where the coefficients of the cohomology group are in $\ZZ_d$.
Applying this procedure to $X=\bar N(\ZZ_d,\Sigma)$ and $[\gamma]=\gamma_b$ produces a linear system that is closely related to the original linear system $(A,b)$. The solution groups of these two linear systems turn out to be isomorphic 
$$
\Gamma(A,b) \xrightarrow{\cong} \Gamma(A_X,b_\gamma),
$$
 thus their solution sets in a group $G$ are in bijective correspondence (Proposition \ref{pro:reduction}).

Motivated by this result, 
we focus on linear systems $(A_X,b_\gamma)$ that come from simplicial sets.
We use the theory of twisted products \cite{May67}. 
Given $(X,\gamma)$ the twisted product is a simplicial set $X_\gamma$ which fits into a fibration sequence
$$
N\ZZ_d\to X_\gamma \to X.
$$ 
The set of $n$-simplices of $X_\gamma$ is given by the product $N\ZZ_d\times X$, but its simplicial structure maps are twisted by $\gamma$.

\begin{customthm}{\ref{thm:pi1-commutative-Xgamma}}
For a (connected) simplicial set $X$, there is a canonical isomorphism of groups
$$
\pi_1(\ZZ_d,X_\gamma) \xrightarrow{\cong} \Gamma(A_X,b_\gamma)
$$
\end{customthm}

The group $\pi_1(\ZZ_d,\cdot)$ is a version of the fundamental group.
Briefly, its generators are the $1$-simplices (edges) of the simplicial set, and simplicial relations come from the $2$-simplices together with additional relations imposing commutativity for the edges on the boundary
and $d$-torsion condition for each generator. This algebraic invariant is first introduced in \cite{okay2020commutative}.
 
{We can associate a linear system to a pair $(G,J_G)$ in a canonical way. Such linear systems are of} particular interest in this paper.
Our starting point is the   central extension
$$
1\to \Span{J_G} \to G \to \bar G\to 1
$$
where $\bar G$ denotes the quotient group.
This extension is classified by a cohomology class $\gamma_G$ in the second cohomology $H^2(\bar G)$ of the quotient group with coefficients in $\ZZ_d$.
This cohomology class can be represented by a $2$-cocycle constructed using a set-theoretic section $\phi:\bar G\to G$ of the quotient homomorphism.
On the level of simplicial sets, the extension above gives rise to a fibration sequence
$$
N\ZZ_d \to N(\ZZ_d,G) \to \bar N(\ZZ_d,G)
$$
classified by a cohomology class similar to $\gamma_G$. The section $\phi$ can be used to construct a representing cocycle {$\gamma_{\phi,d}$ for this class}.
The linear system associated to $(G,J_G)$ is the one associated to the pair $(\bar N(\ZZ_d,G),\gamma_{\phi,d})$.
Writing $(A_G,b_\phi)$ for this linear system we have a group homomorphism 
$$
\Gamma(A_G,b_\phi) \to G.
$$ 
We describe the kernel of this homomorphism using homotopical methods (Proposition \ref{pro:pi1-EZdG}). 

When $d>1$ is an odd integer, in the known cases a linear system admitting a solution in a group $G$ admits a solution in $\ZZ_d$.
We state this as a conjecture (Conjecture \ref{conj:odd}) and provide evidence when $d$ is an odd prime.
 The classes of groups satisfying this conjecture are certain types of extraspecial $p$-groups (Theorem \ref{thm:extraspecial}) and a class of  groups introduced in \cite{frembs2022no} (Theorem \ref{thm:Markus-paper}). Other interesting results we obtain using simplicial methods are as follows:  
\begin{itemize}
\item We show that if a linear system does not admit a solution in $\ZZ_d$ then any simplicial distribution induced by a {density operator}
is contextual (Proposition \ref{pro:linear-system-contextuality}). {The theory of simplicial distributions introduced in \cite{okay2022simplicial} is a framework based on simplicial sets for studying contextuality.}

\item We provide four equivalent characterizations of solutions of linear systems based on group-theoretic and cohomological {criteria}
(Corollary \ref{cor:solution-characterization}). 

\item We consider power maps $\omega_m$ acting on $N(\ZZ_d,G)$ by raising a tuple of group elements to the $m$-th power. Corollary \ref{cor:opposite-odd} applies this map to prove that for odd $d$ certain relations in the solution group simplify. This result appears as the main theorem of \cite{qassim2020classical}.

\item We compute the solution group of the $K_{3,3}$ linear system using our simplicial approach and reproduce some of the key results of \cite{arkhipov2012extending}. 
\end{itemize} 
 
\paragraph{Acknowledgments.}
This work is supported by the US Air Force Office of Scientific Research under award
number FA9550-21-1-0002. The second author would like to thank the Institute for Quantum Computing
for their hospitality during a visit in January 2023, and William Slofstra for fruitful discussions and {providing} the proof of Proposition \ref{pro:unitary-vs-finite-group}.
 
\section{Linear systems}\label{sec:linear-systems}

In this section we will introduce linear systems and their simplicial realizations.
Let $d>1$ be an integer and $\ZZ_d$ denote the additive group of integers modulo $d$.
A {\it linear system} $(A,b)$ over $\ZZ_d$ is a set of linear equations $Ax=b$ specified by a matrix $A\in \ZZ_d^{r}\times \ZZ_d^c$ and a {column} vector $b\in \ZZ_d^c$.
To a linear system we can associate the following data: 
\begin{itemize}
\item  A simplicial complex $\Sigma_A$ with vertex set 
$(\Sigma_A)_0=\set{v_1,v_2,\cdots,v_c}$
and 
maximal simplices
$$
\sigma_i = \set{v_j:\, A_{ij}\neq 0,\; 1\leq j\leq c}
$$
where $1\leq i\leq r$.
\item A collection of functions 
$$A_i:(\Sigma_A)_0 \to \ZZ_d$$
 defined by
$A_i(v_j) = A_{ij}$, where $1\leq i\leq r$.
\item A function 
$$b:\set{\sigma_i:1\leq i\leq r}\to \ZZ_d$$ defined by $b(\sigma_i) = b_i$.
\end{itemize}

\Def{\label{def:solution-in-G}
Let $G$ be a group with a central element $J_G$ of order $d$, i.e., $J_G^d=1$.
A {\it solution} of a linear system $(A,b)$ in the group $G$ is a function $T:(\Sigma_A)_0\to G$ satisfying the following conditions:
\begin{itemize}
\item  $T(v)$ is $d$-torsion, i.e.,  $T(v)^d=1_G$, for all $v\in \Sigma_0$,
\item $\set{T(v):\, v\in \sigma_i}$ pairwise commute for all $1\leq i\leq r$,
\item for all $1\leq i\leq r$ we have
$$
\prod_{v_j\in \sigma_i} T(v_j)^{A_{ij}} = J_G^{b_i}.  
$$ 
\end{itemize} 
We will write $\Sol(A,b;G)$ for the set of solutions of a linear system $(A,b)$ in the group $G$. 
}

Solutions in $\ZZ_d$ in this sense coincide with the ordinary solutions of the linear system in $\ZZ_d$.

For a set $U$ {we will} write $P(U)$ for the power set, that is, the collection of subsets of $U$. 
We will write $\ZZ_d^U$ for the set of functions $U\to \ZZ_d$. This set of functions has a group structure induced by $\ZZ_d$. Note that the rows of $A$ when regarded as functions $A_i:\Sigma_0\to \ZZ_d$ give elements in the group $\ZZ_d^{\Sigma_0}$. We will write $\Span{A_i}$ for the subgroup generated by this function.

\Ex{\label{ex:nonzero-entries}
{\rm
  Let $A$ be a matrix such that the simplicial complex $\Sigma_A$ consists of a unique maximal simplex. That is, $\Sigma_A$ is given by the power set $P(\Sigma_0)$. In other words, $A_{ij}\neq 0$ for all $1\leq i\leq r$ and $1\leq j\leq c$. In this case a solution specified by a function $T:\Sigma_0\to G$ satisfies that $[T(v),T(w)]={T(v)^{-1}T(w)^{-1}T(v)T(w)}=1$ for all vertices $v,w$. We can convert the linear system $(A,b)$ into a row echelon form $(A',b')$ by row operations:
$$
A' = \left( \begin{matrix}
A'' \\ \zero
\end{matrix} \right),\;\;\;\; b' = \left( \begin{matrix}
b_1 \\ b_2
\end{matrix} \right)
$$
where $A''$ contains no rows with all entries zero. Then there is a bijection between $\Sol(A,b;G)$ and $\Sol(A',b';G)$. The latter set is non-empty if $b_2\neq \zero$. Thus
{when studying solutions of linear systems},
 it is a simplifying assumption to require that $\Span{A_i}\neq \Span{A_j}$ for every $i\neq j$.
}
} 

In this paper we will consider linear systems satisfying the following two conditions: 
\begin{enumerate}
\item Each row $A_i$ satisfies $\Span{A_i}\cong \ZZ_d$.
\item For two distinct rows $A_i$ and $A_j$ we have $\Span{A_i}\neq \Span{A_j}$.
\end{enumerate}
These conditions simplify the description of simplicial realizations of linear systems given in Section \ref{sec:SimpReal}. The last property is motivated by Example \ref{ex:nonzero-entries}.
The first property is satisfied by the linear systems of interest, mainly those that come from simplicial sets introduced in Section \ref{sec:LinearSimplicial}.
{Later in Proposition \ref{pro:reduction}, we will show that the solution group of any linear system can be described as the solution group of a linear system satisfying this property.}

\Def{\label{def:solution-group}
{\rm
The {\it solution group} $\Gamma(A,b)$ of a linear system $(A,b)$ is the finitely presented group  generated by $e_v$, where $v\in \Sigma_0$, and $J$ subject to the following relations:
\begin{itemize}
\item $d$-torsion relations: $J^d=e_v^d=1$ for all $v\in \Sigma_0$,
\item commutativity relations: $\set{J,e_v:\,v\in \sigma}$ pairwise commutes for all $\sigma \in \Sigma$,
\item product relations: for all $\sigma_i$ we have
\begin{equation}\label{eq:product-relation}
\prod_{v_j\in \sigma_i} e_{v_j}^{A_{ij}} =J^{b_i} .
\end{equation}
\end{itemize}
}}

Let $\catGrp$ denote the category of groups. We will write $\catGrp(G,H)$ for the set of group homomorphisms. 
{We introduce a category by restricting the morphisms in the category of groups which will be  useful in describing solutions of linear systems.}
Let $\catGrp_J$ denote the {following} category:
\begin{itemize}
\item Objects are pairs $(G,J_G)$ where $J_G\in G$ is a central element of order $d$.
\item A morphism $(G,J_G)\to (H,J_H)$ is given by a group homomorphism $f:G\to H$ such that $f(J_G)=J_H$.
\end{itemize}
We will write $\catGrp_J(G,H)$ for the set of morphisms in this category.

\Pro{
\label{pro:sol-group}
For a linear system $(A,b)$ the following properties hold.
\begin{enumerate}[(1)]
\item There is a bijection 
$$
\Sol(A,b;G) \cong \catGrp_J(\Gamma(A,b),G).
$$

\item The set $\Sol(A,b;G)$ of solutions is non-empty for some $G$  
if and only if $J\in \Gamma(A,b)$ has order $d$. 
\end{enumerate}
}
\Proof{
A solution $T:\Sigma_0 \to G$ can be used to define a group homomorphism $\theta_T:\Gamma(A,b)\to G$ by sending $e_v\mapsto T(v)$ and $J\mapsto J_G$. Conversely, given a group homomorphism in $\catGrp_J(\Gamma(A,b),G)$, restricting to the generators specifies a solution. This proves part (1).

For the second part, assume that the solution set is non-empty for some $G$.
Then the group homomorphism associated to the solution implies that $J$ has order $d$ since the element $\theta_T(J)=J_G$ has order $d$.  
Conversely, if the order of $J$ is $d$ then $T:\Sigma_0\to \Gamma(A,b)$ defined by $T(v)=e_v$ is a solution in  $\Gamma(A,b)$.
}

\Cor{\label{cor:abelian-sol-group}
Assume that $\Gamma(A,b)$ is abelian and $J$ has order $d$. Then $(A,b)$ admits a solution in $\ZZ_d$. 
}

Solutions in the unitary group $U(\CC^m)$, where $m\geq 1$, are of particular importance in the study of linear systems. 
When $m\geq 2$ such solutions are usually referred to as {\it operator solutions}, and when $m=1$ they are called {\it classical solutions}. 
Next, we will show that instead of studying solutions in unitary groups, we can focus on finite groups.

\Pro{[\!\cite{slofstra-private}]  \label{pro:unitary-vs-finite-group}
For a linear system $(A,b)$ we have 
$\Sol(A,b;U(\CC^m))\neq \emptyset$ for some $m\geq 1$ if and only if $\Sol(A,b;G)\neq \emptyset$ for some finite group $G$.
}
\Proof{For a finite group $G$, we can construct a unitary representation of $G$ by inducing the $1$-dimensional representation  
$\Span{J}\xrightarrow{} U(\CC)$ obtained by sending $J$ to $\omega\one$. 
This gives an injective group homomorphism $\phi:G\to U(\CC^m)$ where $m=|G/\Span{J}|$. Then given a solution $T:\Sigma_0\to G$ the composite $\phi\circ T:\Sigma_0\to U(\CC^m)$ is a solution in a unitary group.  

The converse implication follows from the following fact: Any finitely generated subgroup of {the general linear group} $\Gl(\CC^m)$ is residually finite \cite[Theorem 7.116]{dructu2018geometric}.
Therefore the subgroup $G \subset U(\CC^m)$ generated by $\set{T(v_j):\,j=1,\cdots,c}$ is residually finite. Then $G$ is the inverse limit of a sequence of surjective group homomorphisms 
$$\cdots \to G_{i+1} \xrightarrow{f_i} G_i \to \cdots \to G_1$$ 
where each $G_i$ is finite. The element $J\in G$ is represented by a tuple $(J_i)_{i\geq 1}$ of central elements where $J_i\in G_i$. There exists $N \geq 1$ such that $J_N$ has order $d$. Then using the projection $\pi_N:G\to G_N$ we obtain a solution in a finite group. 
}

For a simplicial complex $\Sigma$ we will write $\hat \Sigma$ for the {\it dual complex} consisting of the vertex set 
$$\hat \Sigma_0=\set{\sigma_i:\,1\leq i\leq r}$$ and maximal simplices 
$$\hat \sigma_j = \set{\sigma_i:\, v_j\in \sigma_i},$$
where $1\leq j\leq c$.
An important source of linear systems comes from the incidence matrices of graphs. 
Let $K$ be a graph with vertex set $K_0$ and edge set $K_1$, which we can think of as a simplicial complex.
Let $b:K_0\to \ZZ_d$ be a function. 
Then the incidence matrix $A(K)$, with entries 
$$
A(K)_{v,x}= \left\lbrace
\begin{array}{ll}
1 & v \in x \\
0 & \text{otherwise,}
\end{array}
\right.
$$ 
together with the function $b$ specifies a linear system. Note that $\Sigma_{A(K)}$ coincides with the dual complex $\hat K$. 
{For more on these kinds of linear systems see \cite{paddock2022arkhipov}.}

\begin{figure}[h!]
\centering
\begin{subfigure}{.49\textwidth}
  \centering
  \includegraphics[width=.5\linewidth]{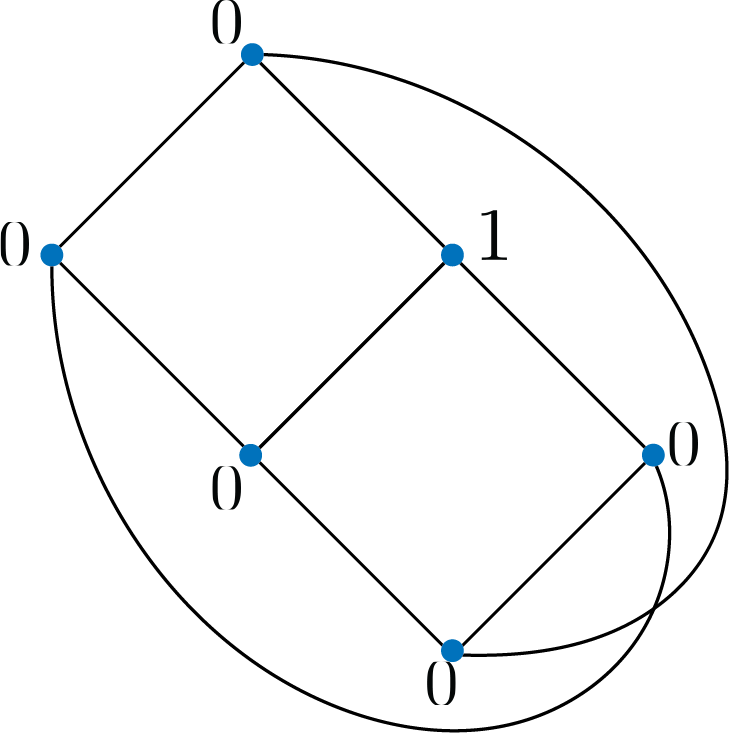}
  \caption{}
  \label{fig:K33-a}
\end{subfigure}%
\begin{subfigure}{.49\textwidth}
  \centering
  \includegraphics[width=.4\linewidth]{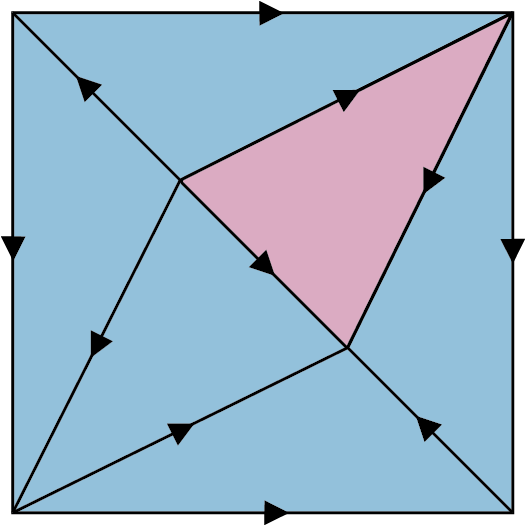}
  \caption{}
  \label{fig:K33-b}
\end{subfigure}
\caption{(a) The $K_{3,3}$ graph. Each vertex is assigned a value $b_i$. (b) Dual of the graph {representing a torus}. Each triangle is assigned a value $b_i$ (pink color corresponds to $1$ value). The top edge is identified with the bottom edge, and the leftmost edge is identified with the rightmost edge.
}
\label{fig:want}
\end{figure}

\Ex{\label{ex:K33}
{\rm
Let $K_{3,3}$ denote the complete bipartite graph illustrated in Figure (\ref{fig:K33-a}). The dual complex $\Sigma_{3,3}=\hat K_{3,3}$ is given by a torus triangulated as in Figure (\ref{fig:K33-b}). It is well known that $K_{3,3}$ admits a solution in $\ZZ_2$ if and only if $\sum_{i=1}^6 b_i=0$; see  \cite{arkhipov2012extending,coladangelo2017robust,paddock2022arkhipov}. In  Section \ref{sec:K33-linear-system} we will study this linear system using our simplicial techniques. 
}
}

\subsection{Simplicial realizations} \label{sec:SimpReal}

In Section \ref{sec:linear-systems} we have seen that a linear system can be described by a simplicial complex, a set of functions each supported on a maximal simplex and a function on the set of maximal simplices. In this section we will express this data in a different way using the language of simplicial sets.

A {\it simplicial set} consists of a sequence of sets $X_1,X_2,\cdots,X_n$ together with 
\begin{itemize}
\item face maps $d_i:X_n\to X_{n-1}$ where $0\leq i\leq n$, and
\item degeneracy maps $s_j:X_n\to X_{n+1}$ where $0\leq j\leq n$
\end{itemize}
satisfying the simplicial identities \cite{GJ99,Friedman12}.
Elements of $X_n$ are called $n$-simplices.
A simplex $\sigma\in X_n$ in the image of $s_j:X_n\to X_{n+1}$ for some $j$ is called {\it degenerate}; otherwise, it is called {\it non-degenerate}.
A map $f:X\to Y$ between two simplicial sets consists of a sequence of functions $f_n:X_n\to Y_n$ for every $n\geq 0$ that respects the face and the degeneracy maps. A simplicial subset of $X$ is a simplicial set $Z$ together with a map $i:Z\to X$ of simplicial sets such that each $i_n$ is given by inclusion of sets $Z_n\subset X_n$.
{
A simplicial set $X$ is {\it connected} if for any two $0$-simplices there is a sequence of $1$-simplices connecting them. That is, for $v,w\in X_0$ there exists $\tau_1,\cdots,\tau_k\in X_1$ such that 
$$
v=d_{i_1}\tau_1,\;\; d_{i_1'}\tau_1 = d_{i_2}\tau_2,\;\; d_{i_2'}\tau_2 = d_{i_3}\tau_2,\;\cdots\;, d_{i'_{k-1}}\tau_{k-1}=  d_{i_k}\tau_k ,\;\; d_{i'k}\tau_k=w
$$ 
where $i_l\in \set{0,1}$ and $i_l\neq i_l'$ for $l=1,\cdots,k$. In this paper we will be mostly concerned with connected simplicial sets.
}

Our examples of simplicial sets will follow a basic construction from algebraic topology known as the nerve construction.
Let $G$ be a group. The nerve space $NG$ is the simplicial set whose set of $n$-simplices is $G^n$ and simplicial structure maps are given by
$$
d_i(g_1,g_2,\cdots,g_n) = \left\lbrace
\begin{array}{ll}
(g_2,\cdots,g_n) & i=0 \\
(g_1,\cdots,g_ig_{i-1},\cdots,g_n) & 0<i<n \\
(g_1,\cdots,g_{n-1})  & i=n
\end{array}
\right.
$$
and
$$
s_j(g_1,g_2,\cdots,g_n)  = (g_1,\cdots,g_j,1,g_{j+1},\cdots,g_n).
$$
For simplicial realizations of linear systems we will be interested in certain simplicial subsets of nerve spaces. 

Recall that for a set $U$ we write $\ZZ_d^U$ for the set of functions $U\to \ZZ_d$, and this set comes with a group structure inherited from $\ZZ_d$. We consider the nerve space $N(\ZZ_d^U)$. Our simplicial realization will be a certain subspace of this nerve space.
Let $\Sigma$ be a simplicial complex. 

\Def{\label{def:simplicial-realization}
{\rm
We define a simplicial
 subset $N(\ZZ_d,\Sigma)  \subset N(\ZZ_d^{\Sigma_0}) $  whose $n$-simplices are given by 
$$
N(\ZZ_d,\Sigma)_n = \set{(s_1,s_2,\cdots,s_n) \in (\ZZ_d^{\Sigma_0})^n  :\, \cup_{i=1}^n\supp(s_i)\in \Sigma}
$$
where $\supp(s)=\set{v\in \Sigma_0:\, s(v)\neq 0}$ denotes  the support of a function $s:\Sigma_0\to \ZZ_d$. 
}}

For a linear system $(A,b)$ with associated simplicial complex $\Sigma$, the  functions $A_i$, where $1\leq i\leq r$, can be regarded as $1$-simplices of $N(\ZZ_d,\Sigma)$. 
Thus $N(\ZZ_d,\Sigma)$ can be thought of as a space that encodes the matrix $A$. Next we introduce another space that can be used for describing solutions of the linear system.
Our construction is a modified version of the nerve construction introduced in \cite{okay2021commutative}; see also \cite{ACT12} for other versions.  

\Def{\label{def:commutative-nerve-d}
{\rm
The {\it mod-$d$ commutative nerve space} is the simplicial subset  $N(\ZZ_d,G)\subset N(G)$ 
whose $n$-simplices are given by 
$$
N(\ZZ_d,G)_n = \set{(g_1,g_2,\cdots,g_n)\in G^n:\, g_i^d=1_G,\; g_ig_j=g_jg_i,\;\forall 1\leq i,j\leq n}.
$$ 
}}

When the $d$-torsion condition is dropped the resulting simplicial set is called the {\it commutative nerve space} and denoted by $N(\ZZ,G)$, a space first introduced in \cite{ACT12}. Its $n$-simplices are given by $n$-tuples of pairwise commuting group elements.

\Lem{\label{lem:maps-into-NG}
A simplicial set map $f:X\to NG$ is determined by its restriction to $1$-simplices, that is, by the function $f_1:X_1\to G$. Moreover, a function $h:X_1\to G$ extends to a simplicial set map $f:X\to NG$ where $f_1=h$ if $h(d_1\sigma)=h(d_2\sigma)h(d_0\sigma)$ for all $\sigma\in X_2$. 
}
\Proof{
See \cite[Proposition 3.13]{okay2022simplicial}.
}

This result can be used to describe maps into a simplicial subspace of $NG$, for example, the mod-$d$ commutative nerve space $N(\ZZ_d,G)$.
The nerve spaces introduced so far are partial monoids in the sense of \cite{broto2021extension}. 
In fact, the theory {introduced there} provides a nice framework for studying the simplicial sets {of interest} in this paper. This approach will be pursued elsewhere. 


{We will need the following simplicial set maps:
\begin{itemize}
\item $\iota:N\ZZ_d \to N(\ZZ_d,G)$ induced by the homomorphism $\ZZ_d \xrightarrow{1\mapsto J_G} G$,
\item $\alpha: \vee_{i=1}^r N\ZZ_d \to N(\ZZ_d,\Sigma)$ whose $i$-th factor is induced by the homomorphism $\ZZ_d \xrightarrow{1\mapsto A_i} \ZZ_d^{\Sigma_0}$, 
\item $\beta: \vee_{i=1}^r N\ZZ_d \to N\ZZ_d$ whose $i$-th factor is induced by the homomorphism $\ZZ_d \xrightarrow{1\mapsto b_i} \ZZ_d$. 
\end{itemize}
}

\begin{pro}\label{pro:solution}
There is a bijective correspondence between  
$\Sol(A,b;G)$
and the set of simplicial set maps $f:N(\ZZ_d,\Sigma)\to N(\ZZ_d,G)$ that make the following diagram commute
\begin{equation}\label{dia:solution}
\begin{tikzcd}[column sep=huge,row sep=large]
\vee_{i=1}^r N\ZZ_d \arrow[r,"\beta"] \arrow[d,"\alpha"',hook]  & N\ZZ_d \arrow[d,"\iota",hook] \\
N(\ZZ_d,\Sigma) \arrow[r,"f"] & N(\ZZ_d,G)
\end{tikzcd}
\end{equation} 
\end{pro} 
\Proof{
We will use Lemma \ref{lem:maps-into-NG}.
A solution $T:\Sigma_0\to G$ can be used to construct a simplicial set map $f:N(\ZZ_d,\Sigma)\to N(\ZZ_d,G)$ by setting $f_1(s)= \prod_{v\in \supp(s)} T(v)^{s(v)}$.  
Conversely, from a simplicial set map we obtain a solution by setting $T(v)=f_1(\delta^v)$.  
}

\subsection{Simplicial distributions}
\label{sec:simplicial-distributions}

Solutions of linear systems play an important role in the study of contextuality \cite{okay2020homotopical}. 
In this section, we will define this notion using the theory of simplicial distributions introduced in \cite{okay2022simplicial} {and consider the simplicial distributions associated with solutions of linear systems.}
 
A distribution on a set $U$ is defined to be a function $p:U\to \RR_{\geq 0}$ with finite support such that $\sum_{u\in U} p(u)=1$.
We write $D(U)$ for the set of distributions on $U$.
There is a function $\delta: U\to D(U)$ sending $u$ to the delta distribution 
$$
\delta^u(u') = 
\left\lbrace
\begin{array}{ll}
1 & u=u'\\
0 & \text{otherwise.}
\end{array}
\right.
$$ 
Given a function $f:U\to V$ {we define} a function $Df: D(U)\to D(V)$ by
$$
Df(p)(v) = \sum_{u\in f^{-1}(v)} p(u). 
$$
For a simplicial set $Y$ the space of distributions  on this simplicial set is defined to be the simplicial set $D(Y)$ with $n$-simplices given by $D(Y_n)$ together with the face and the degeneracy maps
$$
D(d_i): D(Y_n)\to D(Y_{n-1}) \;\;\text{ and }\;\;
D(s_j): D(Y_{n}) \to D(Y_{n+1}).
$$

\Def{\label{def:simp-dist}
A {\it simplicial distribution} is a simplicial set map $p: X\to D(Y)$. We will write $\sDist(X,Y)$ for the set of simplicial distributions.
Given a simplicial set map $s:X\to Y$ we can define a simplicial distribution $\delta^s:X\xrightarrow{s} Y \xrightarrow{\delta} DY$. Simplicial distributions of the form $\delta^s$ are called {\it deterministic distributions}.
We will write $\dDist(X,Y)$ for the set of deterministic distributions.
}

Given a simplicial distribution $p:X\to DY$, we will write $p_\sigma$ for the distribution $p_n(\sigma)\in D(Y_n)$ where {$\sigma\in X_n$}.
We will denote by $\delta:Y\to DY$ the simplicial set map defined by $\delta_n(\sigma)=\delta^{\sigma}$.
There is a canonical map
$$
\Theta: D(\dDist(X,Y)) \to \sDist(X,Y)
$$ 
that sends $d=\sum_{s} \lambda(s) \delta^s$ to the simplicial distribution $\Theta(d)$ defined as follows:
$$
\Theta(d)_\sigma(\theta) = \sum_{r:\,r_\sigma=\theta} \lambda(r)
$$
where $\sigma\in X_n$, $\theta\in Y_n$, and  the summation runs over simplicial set maps $r:X\to Y$ such that $r_\sigma=\theta$.
 
\Def{\label{def:contextuality}
A simplicial distribution $p:X\to DY$ is called {\it contextual} if it does not lie in the image of $\Theta$; otherwise, it is called {\it non-contextual}.
}

A conventional way of formulating contextuality is to use presheaves of distributions \cite{abramsky,abramsky2015contextuality}. Any presheaf of distributions can be realized as a simplicial distribution and the notion of contextuality given in Definition \ref{def:contextuality} specializes to the usual notion formulated in this language \cite{kharoof2022simplicial}.

Contextual simplicial distributions arise from quantum measurements. Let $\hH$ denote a finite dimensional complex Hilbert space. Let $\Proj(
\hH)$ denote the set of projector operators acting on $\hH$, i.e., Hermitian operators that square to itself. A projective measurement is a function $\Pi: U\to \Proj(\hH)$ with finite support such that $\sum_{u\in U} \Pi(u)=\one_\hH$. We write $P_\hH U$ for the set of projective measurements on $U$.
Given a function $f:U\to V$ we define 
$P_\hH f: P_\hH U \to P_\hH V$  by 
$$
P_\hH f(p)(v) = \sum_{u\in f^{-1}(v)} \Pi(u). 
$$ 
For
a simplicial set $Y$ the space of projective measurements 
is given by the simplicial set $P_\hH Y$ whose $n$-simplices are given by $P_\hH Y_n$ together with the simplicial structure maps
$$
P_\hH(d_i): P_\hH(Y_n)\to P_\hH(Y_{n-1}) \;\;\text{ and }\;\;
P_\hH(s_j): P_\hH(Y_{n}) \to P_\hH(Y_{n+1}).
$$ 
Given a density operator $\rho$, a positive operator of trace $1$, we can define a simplicial set map
$
\rho_*: P_\hH Y\to DY
$
by sending $\Pi\in P_\hH Y_n$ to the distribution 
$$
\rho_*\Pi (\theta) = \Tr(\rho \Pi(\theta)). 
$$
Checking that this is indeed a simplicial set map is straight-forward, and essentially follows from the linearity of the trace.

\Pro{\label{pro:spectral-dec} 
The spectral decomposition theorem gives an isomorphism of simplicial sets
$$
\sd:N(\ZZ_d,U(\hH)) \to P_\hH N\ZZ_d.
$$
}
\Proof{
 In degree $n$, $\sd$ is described as follows: A  tuple $(A_1,\cdots,A_n)$ of unitary matrices to the projective measurement $\Pi:\ZZ_d^n\to \Proj(\hH)$ where $\Pi(a_1,\cdots,a_n)$ projects onto the simultaneous eigenspace with eigenvalues $(\omega^{a_1},\cdots,\omega^{a_n})$.
For details see \cite[Proposition 6.3]{okay2022simplicial}.
}

Now, let $G$ be a group and $\chi:G\to U(\hH)$ be a group homomorphism, in other words, a unitary representation.
We have a commutative diagram of simplicial sets
$$
\begin{tikzcd}[column sep=huge,row sep=large]
N(\ZZ_d,G) \arrow[r,"\chi_*"] \arrow[dr,"p_\rho"'] &  P_\hH N\ZZ_d \arrow[d,"\rho_*"] \\
& D(N\ZZ_d) 
\end{tikzcd}
$$
where $\chi_*$ in degree $n$ sends $(g_1,\cdots,g_n)$ to $\sd(\chi(g_1),\cdots,\chi(g_n))$.

\Pro{\label{pro:linear-system-contextuality}
Let $(A,b)$ be a linear system over $\ZZ_d$ and $\Sigma$ the associated simplicial complex. Let $\chi:G\to U(\hH)$ be a group homomorphism.
Assume that $(A,b)$ admits a solution in $G$, but not in $\ZZ_d$. Let $f:N(\ZZ_d,\Sigma)\to N(\ZZ_d,G)$ denote the simplicial set map corresponding to the solution (Proposition \ref{pro:solution}). 
Then the simplicial distribution given by the composite
$$
N(\ZZ_d,\Sigma) \xrightarrow{f} N(\ZZ_d,G) \xrightarrow{p_\rho} DN\ZZ_d 
$$
is contextual for any density operator $\rho$.
} 
\begin{proof}
In the contrary assume that $p=f\circ p_\rho$ is non-contextual, that is, there exists $d=\sum_r \lambda(r) \delta^r$ in $D(\dDist(X,N\ZZ_d))$, where $X=N(\ZZ_d,\Sigma)$, such that $\Theta(d)=p$. 
Let $s:X\to N\ZZ_d$ be such that $\lambda(s)\neq 0$. Then for all $\sigma\in X_n$ we have
\begin{equation}\label{eq:supp-s}
p_\sigma(s_\sigma) = \sum_{r:\, r_\sigma=s_\sigma} \lambda(r)   > \lambda(s) >0.
\end{equation}
Consider the commutative diagram
$$
\begin{tikzcd}[column sep=huge,row sep=large]
\vee_{i=1}^r N\ZZ_d \arrow[r,"\beta"] \arrow[d,"\alpha"',hook]  & N\ZZ_d \arrow[d,"1\mapsto J_G",hook] \arrow[r,equal] & N\ZZ_d  \arrow[d,"\delta",hook] \\
N(\ZZ_d,\Sigma) \arrow[r,"f"]  \arrow[ru,"s",dashed] & N(\ZZ_d,G) \arrow[r,"p_\rho"] & DN\ZZ_d
\end{tikzcd}
$$
The left square is Diagram \ref{dia:solution}.
Our goal is to show that $s$ makes the left-top triangle in the left square commute, i.e., $s\circ\alpha=\beta$. By Lemma \ref{lem:maps-into-NG} it suffices to verify commutativity in degree $1$. Let $1\in \ZZ_d$ denote the $1$-simplex in the $i$th factor of the wedge product $\vee^r N\ZZ_d$. 
By the vertical map this element maps to the $1$-simplex $A_i$ and under the horizontal map it maps to $b_i$. By commutativity of the outer square we obtain that $p_{A_i}=\delta^{b_i}$.
Then by Equation (\ref{eq:supp-s}) we have
$$
\delta^{b_i}(s_{A_i})= p_{A_i}(s_{A_i})>0. 
$$
That is, $s_{A_i}=b_i$. Therefore $s\circ\alpha=\beta$.
\end{proof}

This result explains the importance of operator solutions of linear systems. Those which admits solution in a group but not in $\ZZ_d$, i.e., solutions in the conventional sense, give rise to contextual simplicial distributions for any density operator (quantum state). Our approach is based on simplicial sets to make the connection to the theory of simplicial distributions in a more direct way. One can pass through presehaves of distributions to specialize the corresponding result in that language. That is, such linear systems also give rise to contextual presheaves of distributions \cite{abramsky}.

\subsection{Linear systems from simplicial sets} \label{sec:LinearSimplicial}
 
Our goal in this section is to associate a linear system to a given simplicial set. We will show that as far as the solution groups are concerned  any linear system can be converted to one that comes from a simplicial set.

{Next we introduce cohomology of simplicial sets. Throughout the paper we restrict to cohomology with coefficients in $\ZZ_d$.}
Given a simplicial set an $n$-cochain taking values in $\ZZ_d$ is a function $X_n\to \ZZ_d$. We will write $C^n(X)$ for the set of $n$-cochains. There is a coboundary map
$$
d_n: C^{n-1}(X) \to C^{n}
(X)
$$
defined by sending $f:X_{n-1}\to \ZZ_d$ to he function $d_nf:X_{n}\to \ZZ_d$:
$$
d_nf(\sigma) = \sum_{i=0}^n f(d_i\sigma)
$$
for $\sigma\in X_{n}$. The $n$-th cohomology group is defined by the quotient group
$$
H^n(X) = \frac{\ker (d_{n+1})}{\im(d_n)}.
$$
{Sometimes we will write $H^n(X,\ZZ_d)$ to emphasize the coefficients.}
We will construct a linear system  associated to 
a simplicial set together with a $2$-cochain.

\Def{\label{def:linear-system-of-simplicial-set}
Given a simplicial set $X$ and a $2$-cochain $\gamma:X_2\to \ZZ_d$ we define a linear system $(A_X,b_\gamma)$ where $A$ is a $|X_2|\times |X_1|$ matrix and $b$ is a column vector of size $|X_2|$:
$$
A_{\sigma,x} = 
\left\lbrace
\begin{array}{ll}
\sum_{d_i\sigma=x} (-1)^i & x\in \partial\sigma\\ 
0 & \text{otherwise,}
\end{array}
\right.
\;\;\text{ and } \;\;
b_\sigma = -\gamma(\sigma)
$$
where the summation runs over $0\leq i\leq 3$ such that $d_i\sigma=x$. 
}

Our main examples in this paper will come from simplicial sets satisfying the property that $|\partial\sigma|=3$ for every $2$-simplex of $X$. In this case the definition of $A$ becomes:
$$
A_{\sigma,x} =  
\left\lbrace
\begin{array}{ll}
1 &  x\in \set{d_0\sigma,d_2\sigma} \\
-1 & x=d_1\sigma \\
0 & \text{otherwise.}
\end{array}
\right.
$$
A solution in $G$ for the  linear system $(A_X,b_\gamma)$ consists of a function $T:X_1\to G$ that satisfies 
\begin{itemize}
\item $T(x)^d=1$ for all $x\in X_1$,
\item $\set{T(d_i\sigma):\,i=0,1,2}$ pairwise commute for all $\sigma\in X_2$,
\item for every $\sigma\in X_2$ we have
$$
T(d_2\sigma) T(d_0\sigma)T(d_1\sigma)^{-1} = J_G^{-\gamma(\sigma)} .
$$
\end{itemize}

The simplicial set 
$N(\ZZ_d,\Sigma)$ (Definition \ref{def:simplicial-realization}) 
{is our main example.}
As we have seen each row $A_i$ can be regarded as an element in $\ZZ_d^{\Sigma_0}$. The nerve space $N\Span{A_i}$ is contained in $N(\ZZ_d,\Sigma)$. The union of these simplicial subsets over $1\leq i\leq r$ is the simplicial subset given by the wedge $\vee_{i=1}^r N\Span{A_i}$. 
Here we use the wedge notation to emphasize that for $i\neq j$ the intersection of $N\Span{A_i}\cap N\Span{A_j}$ is given by {$\Delta^0$, the simplicial set with $(\Delta^0)_n=\set{\ast}$ for $n\geq 0$ representing a point.}


\Def{\label{def:Nbar}
Let $\bar N(\ZZ_d,\Sigma)$ denote the simplicial set obtained by taking the quotient of $N(\ZZ_d,\Sigma)$ by the subspace $\vee_{i=1}^r N\Span{A_i}$. More explicitly, the  $n$-simplices are given by
$$
\bar N(\ZZ_d,\Sigma)_n =  \set{\bar 0}\, \sqcup\, 
\left( N(\ZZ_d,\Sigma)_n -  \vee_{i=1}^r (N\Span{A_i})_n \right).
$$
}

Next we construct a cohomology class  
$\gamma_b$. 
We will use the cohomology long exact sequence of the cofiber sequence 
$$
\vee_{i=1}^r N\ZZ_d  \to N(\ZZ_d,\Sigma) \to \bar N(\ZZ_d,\Sigma)
$$
where we identify $\Span{A_i}$ with $\ZZ_d$.
There is an associated long exact sequence in mod $d$ cohomology
\begin{equation}\label{eq:connecting-hom}
\cdots \to H^1(\bar N(\ZZ_d,\Sigma)) \to H^1(N(\ZZ_d,\Sigma)) \to  H^1(\vee_i N\ZZ_d) \xrightarrow{\delta} H^2( \bar N(\ZZ_d,\Sigma) ) \to \cdots
\end{equation}
The vector $b\in \ZZ_d^r$ can be identified with an element in $H^1(\vee_i N\ZZ_d)$ via the following isomorphism 
\begin{equation}\label{eq:H1-N}
H^1(\vee_i N\ZZ_d) \cong  \ZZ_d^r.
\end{equation} 
Let $\gamma_b$ be the $2$-cocycle on $\bar N(\ZZ_d,\Sigma)$
defined by 
\begin{equation}\label{eq:gamma-b}
\gamma_b(\bar s_1,\bar s_2) = d\tilde b(\bar s_1,\bar s_2)
\end{equation}
where $\tilde b$ is the $1$-cochain on $N(\ZZ_d,\Sigma)$ 
given by
\begin{equation}\label{eq:tilde-b}
\tilde b(s) = \left\lbrace
\begin{array}{cc}
a b_i & s = a A_i \\
0 & \text{otherwise.}
\end{array}
\right.
\end{equation}
In Lemma \ref{lem:delta-b} we show that $\delta(b)=\gamma_b$.

\Pro{\label{pro:reduction}
For a linear system $(A,b)$, let $X=\bar N(\ZZ_d,\Sigma)$ and $\gamma$ be such that $[\gamma]=\gamma_b$. Then there is an isomorphism between the solution groups
$$
\Gamma(A,b) \to \Gamma(A_X,b_\gamma).
$$
}
\Proof{
Proof is given in Section \ref{sec:proof-reduction}.
}

\begin{ex}
{\rm
Consider the following linear system over $\mathbb{Z}_2$:
\[
\begin{pmatrix}
1&1\\
1&0
\end{pmatrix}
\begin{pmatrix}
x\\
y
\end{pmatrix}
=
\begin{pmatrix}
b_1\\
b_2
\end{pmatrix}.
\]
The associated  simplicial complex is given by $\Sigma=P(\Sigma_0)$ where $\Sigma_0=\set{v_1,v_2}$. 
The two rows can be regarded as the functions: $A_1\colon v_1\mapsto 1,v_2\mapsto 1$ and $A_2\colon v_1\mapsto 1,v_2\mapsto 0$.
It is notationally convenient to 
encode functions $\Sigma\to\mathbb{Z}_2$ by two 'bits', representing values on $v_1$ and $v_2$ respectively. Thus we encode $A_1$ as $11$ and $A_2$ as $10$.
The simplicial set $N(\mathbb{Z}_2,\Sigma)$ consists of functions $\{00,01,10,11\}$ as $1$-simplices.  
The $n$-simplices of $N(\mathbb{Z}_2,\Sigma)$ are given by $n$-tuples $(s_1,s_2,\ldots,s_n)\in N(\mathbb{Z}_2,\Sigma)_1^n$ with no further restrictions.

The simplicial set $\bar{N}(\mathbb{Z}_2,\Sigma)$ is obtained by identifying  
the two copies of $N\ZZ_2$ to the base point. 
The set of $1$-simplices of $\bar{N}(\mathbb{Z}_2,\Sigma)$ is given by $\{\bar0=00,01\}$ and the set of $2$-simplices is 
$$\{\bar0=(00,00),(00,01),(01,00),(01,01),(01,10),(10,01),(11,01),(01,11),(11,10),(10,11)\}.$$ 
Note that in $\bar{N}(\mathbb{Z}_2,\Sigma)_2$ only the elements of the form $(aA_i,bA_i)$ for $i=1,2$ and $a,b\in\mathbb{Z}_2$ are identified to the base point.
In this case the linear system $(A_X,b_\gamma)$ where 
$X=\bar{N}(\mathbb{Z}_2,\Sigma)$ takes the following form:
\[
\begin{blockarray}{ccc}
&00&01\\
\begin{block}{c(cc)}
(00,00)&1&0\\
(00,01)&1&0\\
(01,00)&1&0\\
(01,01)&1&0\\
(01,10)&0&1\\
(10,01)&0&1\\
(11,01)&0&1\\
(01,11)&0&1\\
(11,10)&0&1\\
(10,11)&0&1\\
\end{block}
\end{blockarray}
=
\begin{pmatrix}
0\\
0\\
0\\
0\\
b_1+b_2\\
b_1+b_2\\
b_2+b_1\\
b_2+b_1\\
b_1+b_2\\
b_1+b_2\\
\end{pmatrix}
\]
This linear system can be reduced to the initial one.
} 
\end{ex}

\section{Twisted products}\label{sec:TwistedProds}

Let $X$ be a simplicial set and $\gamma:X_2\to \ZZ_d$ be a $2$-cocycle. 
We will assume that 
$\gamma$ is normalized in the sense that 
\begin{equation}\label{eq:normalized}
\gamma(s_iX_1)=0\;\; \text{ for }\;i=0,1.
\end{equation}
To this cocycle we can associate a fibration
$$
N\ZZ_d \xrightarrow{i} X_\gamma \to X.
$$
The total space $X_\gamma$ can be described explicitly using twisted products.  

\Def{\label{def:twisted-prod}
{\rm  
The twisted product $X_\gamma=N\ZZ_d\times_{\gamma} X$ is a simplicial set defined as follows:
\begin{itemize}
\item the set of $n$-simplices is given by $\ZZ_d^n\times X_n$,
\item for $\alpha\in \ZZ_d^n$ and $\tau\in X_n$ the simplicial structure maps are given by
$$
d_i(\alpha,\tau) = (d_i\alpha,d_i\tau),\;\;\;s_j(\alpha,\tau)=(s_j\alpha,s_j\tau)
$$
where $1\leq i\leq n$ and $0\leq j\leq n$; and when $i=0$ we have
$$
d_0(\alpha,\tau) = ( \eta(\tau)+d_0\alpha,d_0\tau )
$$
where $\eta:X_n\to \ZZ_d^{n-1}$ is the twisting function that depends on the cocycle $\gamma$; see Equation (\ref{eq:eta}) {in Section \ref{sec:ClassFib}.}
\end{itemize}
}}

\subsection{Commutative fundamental group}  \label{sec:CommutativeFund}

Let  
 $\catsSet$ denote  
the category of  simplicial sets. 
 There is a well-known adjunction $\pi_1 \adjoint N$ between the nerve functor $N:\catGrp\to \catsSet$ {that sends} a group to the nerve space and the functor $\pi_1:\catsSet\to \catGrp$ {that sends} a simplicial set $X$ to the group defined as follows:
$$
\pi_1(X) = \Span{ e_{x},\,x\in X_1:\, e_{d_2\sigma} e_{d_0\sigma} = e_{d_1\sigma}\;\forall \sigma\in X_2 },
$$  
We will refer to $\pi_1X$ 
as the {\it algebraic fundamental group}. 
Note that this group coincides with the fundamental group of the geometric realization for reduced simplicial sets, 
i.e., $X_0=\set{\ast}$. 
For example, this is the case for $N(\ZZ_d,G)$ and $N(\ZZ_d,\Sigma)$.

\Lem{\label{lem:deg-1}
For every $x\in s_0X_0$ we have $e_x=1$. 
}
\Proof{
Consider the $2$-simplex $\sigma =s_0x$. We have $d_i\sigma=x$ for all $i=0,1,2$ by the simplicial relations. Therefore 
$
e_{x} e_{x} = e_{x},
$
which implies that $e_x=1$.
}

By the adjunction 
there is a natural isomorphism
$$
\catsSet(X,NG)\cong \catGrp(\pi_1 X,G).
$$ 
Recall that $N(\ZZ_d,G)$ is a simplicial subset of $NG$, and there is a similar adjunction that can be used to provide an algebraic description of simplicial set maps. Next, we describe this adjunction. To this end,
we introduce a version of the algebraic fundamental group. 

\Def{\label{def:comm-fund-group}
{\rm
The {\it commutative ($d$-torsion algebraic) fundamental group} of a simplicial set $X$ is the group  $\pi_1(\ZZ_d,X)$ generated by $e_x$ for $x\in X_1$ subject to the relations
\begin{itemize}
\item $e_x^d=1_G$ for all $x\in X_1$,
\item $[e_{d_i\sigma},e_{d_j\sigma}]=1_G$ for all $\sigma\in X_2$ and $i,j=0,1,2$,
\item $e_{d_2\sigma} e_{d_0\sigma} = e_{d_1\sigma}$ for all $\sigma\in X_2$.
\end{itemize}
}}

As a consequence of this definition there is a surjective group homomorphism $\pi_1X \to \pi_1(\ZZ_d,X)$ defined by the identity map on the set of generators.

\Lem{\label{lem:fund-NZdG}
The fundamental group of $X=N(\ZZ_d,G)$ has the following presentation
$$
\Span{e_g,\, g\in G_{(d)}:\, e_ge_h = e_{gh}\text{ whenever } [g,h]=1 }.
$$
In particular,  the quotient homomorphism $\pi_1(X)\to \pi_1(\ZZ_d,X)$ is an isomorphism.
} 
\Proof{
Using $\sigma=(g,h)$, for $g,h\in G_{(d)}$ such that $[g,h]=1$, we obtain the relation $e_g e_h = e_{gh}$.
From which it follows that $e_{g}^d= e_{g^d}=e_1=1$ by Lemma \ref{lem:deg-1} and 
$
e_{g}e_h = e_{gh}=e_{hg} = e_h e_g
$. 
}

\Lem{\label{lem:maps-into-NZdG}
Simplicial set maps $f:X\to N(\ZZ_d,G)$ are in bijective correspondence with  
functions $h:X_1\to G_{(d)}$ 
satisfying {$h(d_1\sigma)=h(d_2\sigma)h(d_0\sigma)$} for all $\sigma\in X_2$.
} 
\Proof{
Follows from Lemma \ref{lem:maps-into-NG}
}

\Pro{\label{pro:adjunction}
The functors $\pi_1(\ZZ_d,\cdot)$ and $N(\ZZ_d,\cdot)$ constitute an adjoint pair, i.e.,
there is a natural bijection 
$$
\catsSet(X,N(\ZZ_d,G))\cong \catGrp(\pi_1(\ZZ_d, X),G).
$$ 
}
\Proof{
A simplicial set map $f:X\to NG$ factors through the inclusion $N(\ZZ_d,G)\subset NG$ if and only if the adjoint map $\hat f: \pi_1 X\to G$ factors through the quotient map $\pi_1 X \to \pi_1(\ZZ_d, X)$. This follows from Lemma \ref{lem:fund-NZdG} and Lemma \ref{lem:maps-into-NZdG}.
}

Also, the  fundamental group of $N(\ZZ_d,\Sigma)$ has a nice description. 
For a simplex $\sigma\in \Sigma$ let us define the function  $\delta^\sigma:\Sigma\to \ZZ_d$ by
$$
\delta^\sigma(v) = \left\lbrace
\begin{array}{ll}
1 & v\in \sigma\\
0 & \text{otherwise.}
\end{array}
\right.
$$
When $\sigma=\set{v}$ we will abbreviate the notation as $\delta^v$ and 
write $e_v$ for $e_{\delta^v}$. 

\Lem{\label{lem:pi1-NZdSigma}
The fundamental group of $X=N(\ZZ_d,\Sigma)$ has the following presentation
$$
\Span{e_v,\, v\in \Sigma_0:\, e_ve_w = e_{\delta^{\set{v,w}} }\text{ whenever } \set{v,w}\in \Sigma }.
$$
Moreover,   the quotient homomorphism $\pi_1(X)\to \pi_1(\ZZ_d,X)$ is an isomorphism.
}
\Proof{
Note that a function $s:\Sigma_0\to \ZZ_d$   can be written as $s= \sum_{v\in \Sigma_0} k_v \delta^v$. Using simplices of the form $\sigma=(k\delta^v,l\delta^w)$, where $v,w$ is such that $\set{v,w}\in \Sigma$, we can write
\begin{equation}\label{eq:e_s}
e_s = \prod_{v\in \Sigma_0} e_{\delta^v}^{k_v}.
\end{equation}
Then each generator $e_{\delta^v}$ is $d$-torsion since $e_{\delta^v}^d = e_{d\delta^v}=e_0=1$. Moreover, $[e_{\delta^v},e_{\delta^w}]=1$ whenever $\set{v,w}\in \Sigma$.
}

\Cor{\label{cor:sol-fund}
The set $\Sol(A,b;G)$ of solutions  
is in bijection  
with the set of group homomorphisms $\theta:\pi_1N(\ZZ_d,\Sigma)\to G$ satisfying $\theta(e_{A_i})=J_G^{b_i}$ for all $1\leq i \leq r$.
}
\begin{proof}
By the adjunction of Proposition \ref{pro:adjunction} and the computation in Lemma \ref{lem:pi1-NZdSigma}, 
the diagrams in Proposition \ref{pro:solution} characterizing solutions are in bijective correspondence with diagrams of group homomorphisms of the form 
$$
\begin{tikzcd}[column sep=huge,row sep=large]
\ast_{i=1}^r \ZZ_d \arrow[d,"\ast_i (1\mapsto e_{A_i})"'] \arrow[r,"\ast_i (1\mapsto b_i)"] & \ZZ_d \arrow[d,"1\mapsto J_G"] \\
\pi_1N(\ZZ_d,\Sigma) \arrow[r,"\theta"] & G 
\end{tikzcd}
$$
\end{proof}

Note that for a group homomorphism $\theta:\pi_1N(\ZZ_d,G)\to G$ such that $\theta(e_{A_i})=J_G^{b_i}$ the images $\theta(e_v)$ will be $d$-torsion and elements in $\set{\theta(e_v):v\in \sigma_i}$ will pairwise commute. Using Equation (\ref{eq:e_s}) the condition  $\theta(e_{A_i})=J_G^{b_i}$ is equivalent to 
$$
\theta(e_{A_i}) = \prod_{v\in \sigma_i} \theta(e_v)^{A_{ij}} = J_G^{b_i}.
$$

\subsection{Fundamental group of twisted products}

{In this section}
we will 
describe the algebraic fundamental group of $X_\gamma$.  
Using the explicit formula for $\eta$ given in Equation (\ref{eq:eta}) in Section \ref{sec:ClassFib} we have
$$
d_0(a,b;\sigma) = (\gamma(\sigma)+ b;d_0\sigma)
$$
where $\sigma\in X_2$ and $(a,b)\in \ZZ_d^2$.
The fundamental group $\pi_1(X_\gamma)$ is generated by $e_{a,x}$, where $(a,x)\in \ZZ_d\times X_1$, subject to the relations
\begin{equation*} 
e_{a+b,d_1\sigma} = e_{a,d_2\sigma}\, e_{\gamma(\sigma)+b,d_0\sigma}
\end{equation*}
where $((a,b),\sigma)\in \ZZ_d^2\times X_2$.  
Here are some useful relations:
\begin{enumerate}[(1)] 
\item For $\sigma=(a,b;s_0x)$, where $x\in X_1$,
\begin{equation*} 
e_{a+b,x} = e_{a,s_0(d_1x)} e_{b,x}
\end{equation*}
where $\gamma(s_0x)=0$ by the normalization condition.
In particular, when $b=0$ this gives 
\begin{equation}\label{eq:s0x-b=0}
e_{a,x}=e_{a,s_0(d_1x)}e_{0,x}
\end{equation}
On the other hand, setting $x=s_0v$ gives 
\begin{equation*}
e_{a+b,s_0v} = e_{a,s_0v} e_{b,s_0v},
\end{equation*}
using which we obtain $e_{a,s_0v}^d=1$ and $[e_{a,s_0v},e_{b,s_0v}]=1$.
\item For $\sigma=(a,b;s_1x)$, where $x\in X_1$,
\begin{equation*}
e_{a+b,x} = e_{a,x} e_{b,s_0(d_0x)}
\end{equation*}
where $\gamma(s_1x)=0$ by the normalization condition.
In particular, when $a=0$ this gives 
\begin{equation}\label{eq:s1x-a=0}
e_{b,x}=e_{0,x}e_{b,s_0(d_0x)}.
\end{equation}
\end{enumerate}

We will write $e_{x}=e_{0,x}$. Combining (1) and (2) we obtain the following.  

\Lem{\label{lem:comm-rel}
We have
\begin{equation}\label{eq:comm-rel}
e_{a,s_0(d_1x)} =  e_{x}\,  e_{a,s_0(d_0x)}\, e_{x}^{-1}.
\end{equation}
}
\Proof{
Follows from Equation (\ref{eq:s0x-b=0}) and (\ref{eq:s1x-a=0}).
}

\begin{figure}[h!] 
\centering
\includegraphics[width=.25\linewidth]{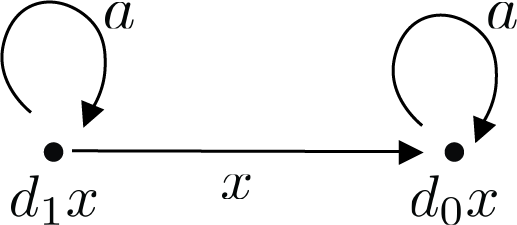}
\caption{The generator $e_{a,s_0(d_1x)}$ corresponds to the loop $a$ at the vertex $d_1x$, the source of the edge $x$. The conjugate element $e_{x}\,  e_{a,s_0(d_0x)}\, e_{x}^{-1}$ is the loop obtained by traversing $x$, $a$, and $x$ in the reverse direction. 
There two elements coincide in $\pi_1X_\gamma$.
}
\label{fig:conjugation}
\end{figure}

See Figure (\ref{fig:conjugation}) for a pictorial representation of this generator.  

\Lem{\label{lem:central}
If $X$ is reduced then $e_a=e_{a,s_0(\ast)}$ is central. 
}
\Proof{
Using  
Equation (\ref{eq:s1x-a=0}) we can write $e_{a,x}=e_{x}e_{a}$. Then the result follows from   
Equation (\ref{eq:s0x-b=0}) and (\ref{eq:comm-rel}).
}

Now, we turn to the commutative fundamental group. The additional relations are 
\begin{enumerate}[(i)]
\item $e_{a,x}^d=1$ for every $(a,x)\in \ZZ_d\times X_1$,
\item $\set{e_{\gamma(\sigma)+b,d_0\sigma} , e_{a+b,d_1\sigma} , e_{a,d_2\sigma}}$ pairwise commute for every $(a,b;\sigma)\in \ZZ_d^2\times X_2$.
\end{enumerate}
When $X$ is connected  we will write 
\begin{equation}\label{eq:e-a}
e_a=e_{a,s_0v}
\end{equation}
 since in the commutative fundamental group the generators $e_{a,s_0v}$ with different $v$'s are all identified.

\Thm{\label{thm:pi1-commutative-Xgamma}
Let $X$ be a connected simplicial set. 
The assignment $e_{a,x}\mapsto J^ae_x$ defines an isomorphism of  groups  
$$
\pi_1(\ZZ_d,X_\gamma) \to \Gamma(A_X,b_{\gamma}).
$$ 
Moreover, for two cocycles $\gamma$ and $\gamma'$ such that $[\gamma]=[\gamma']$ we have an isomorphism $\Gamma(A_X,b_{\gamma})\cong \Gamma(A_X,b_{\gamma'})$.
}
\Proof{ 
We will denote the homomorphism in the statement by $\phi$. We firstly prove that it is a well-defined homomorphism of groups, i.e., the relators in $\pi_1(\zd,X_\gamma)$ are taken to $1$ by $\phi$.
\begin{enumerate}
\item Let $(a,b;\sigma)$ be a $2$-simplex in $X_\gamma$. Consider the relation $e_{a,d_2\sigma}e_{\gamma(\sigma)+b,d_0\sigma}e^{-1}_{a+b,d_1\sigma}=1$.  
We have
\[
\begin{aligned}
\phi\left(e_{a,d_2\sigma}e_{\gamma(\sigma)+b,d_0\sigma}e^{-1}_{a+b,d_1\sigma}\right) &= J^ae_{d_x\sigma}J^{\gamma(\sigma)+b}e_{d_0\sigma}J^{-a-b}e^{-1}_{d_1\sigma}\\
&=J^{\gamma(\sigma)}e_{d_x\sigma}e_{d_0\sigma}e^{-1}_{d_1\sigma} \\
&=J^{\gamma(\sigma)}J^{-\gamma(\sigma)}=1,
\end{aligned}
\]
where we use the  
centrality of the element $J$ and the product relation  
$e_{d_x\sigma}e_{d_0\sigma}e^{-1}_{d_1\sigma}=J^{-\gamma(\sigma)}$ in the group $\Gamma(A_X,b_\gamma)$ (Definition \ref{def:solution-group} and \ref{def:linear-system-of-simplicial-set}).

\item Consider the relation $e_{a,x}^d=1$. We have 
\[
\phi(e_{a,x}^d)=(J^ae_{x})^d=J^{ad}e_{x}^d=1,
\]
where we used  the centrality of the element $J$ and the fact that all generators in $\Gamma(A_X,b_\gamma)$ are of order $d$.

\item The last relation is that for a $2$-simplex $(a,b;\sigma)\in X_\gamma$ the elements $\{e^{-1}_{a+b,d_1\sigma},e_{a,d_2\sigma},e_{\gamma(\sigma)+b,d_0\sigma}\}$ pairwise commute. Note that image of every element from this set is the product of a power of $J$ and an element $e_{d_i\sigma}$ for $i=0,1,2$. However, the elements $e_{d_0\sigma},e_{d_1\sigma}$ and $e_{d_2\sigma}$ pairwise commute in $\Gamma(A_X,b_\gamma$) and  $J$ is central. Thus $\phi(e^{-1}_{a+b,d_1\sigma}),\phi(e_{a,d_2\sigma}),\phi(e_{\gamma(\sigma)+b,d_0\sigma})$ also pairwise commute.
\end{enumerate}
Now, we define  $\psi\colon \Gamma(A_X,b_\gamma)\to \pi_1(\zd,X_\gamma)$ by $\psi(J)=e_{1}$ and $\psi(e_x)=e_x$. Note that if $\psi$ is a well-defined group homomorphism then $\psi$ and $\phi$ are mutually inverse, proving that $\phi$ is an isomorphism. So we need to check again that $\psi$ is taking the relators of $\Gamma(A_X,b_\gamma)$ to $1$ in $\pi_1(\zd,X_\gamma)$.
\begin{enumerate}
\item By definition, the elements $e_x\in\pi(\zd,X_\gamma)$ are $d$-torsion, proving that $\psi(e^d_x)=1$. We also have that the element $e_1=e_{1,s_0(v)}$ for some vertex $v\in\left(X_\gamma\right)_0$ is $d$-torsion, so $\psi(J^d)=1$.
\item By Lemma \ref{lem:central}  
the element $\psi(J)=e_1$ is central in $\pi_1(\zd,X_\gamma)$.

\item Let $\sigma$ be a $2$-simplex in $X$. We need to show that $\psi(e_{d_i\sigma})$ pairwise commute in $\pi(\zd,X_\sigma)$ for $i=0,1,2$. There is a $2$-simplex $(0,0;\sigma)$ in $X_\gamma$, showing that these elements indeed commute.

\item Finally, we need to consider the product relation. Let $\sigma$ be a $2$-simplex in $X$. Then we need to show that
\[
\psi\left(J^{\gamma(\sigma)}e_{d_0\sigma}e_{d_2\sigma}e^{-1}_{d_1\sigma}\right)=1.
\]
Consider the $2$-simplex $(0,0;\sigma)$ in $X_\gamma$. It gives the following relation in $\pi_1(\zd,X_\gamma)$:
\[
e_{\gamma(\sigma),d_0\sigma}e_{d_2\sigma}e_{d_1\sigma}^{-1}=1.
\] 
However, we can write the first element in this relation as $e_{\gamma(\sigma),d_0\sigma}=e_{d_0\sigma}e_{\gamma(\sigma)}$
and by the centrality of the element $e_1$ we obtain that $e_{d_0\sigma}e_{d_2\sigma}e^{-1}_{d_1\sigma}=e^{-\gamma(\sigma)}.$
Combining these gives 
\[
\psi\left(J^{\gamma(\sigma)}e_{d_0\sigma}e_{d_2\sigma}e^{-1}_{d_1\sigma}\right)=e^{\gamma(\sigma)}_1e_{d_0\sigma}e_{d_2\sigma}e^{-1}_{d_1\sigma}=e^{\gamma(\sigma)}e^{-\gamma(\sigma)}=1.
\]
\end{enumerate}
We   conclude that $\psi$ is a well-defined group homomorphism. Thus the   result follows.

For the second statement, observe that Theorem \ref{thm:classification-fib} implies that when the cohomology classes $[\gamma]$ and $[\gamma']$ coincide there is an isomorphism of bundles. Therefore $\pi_1(\ZZ_d,X_\gamma)\cong \pi_1(\ZZ_d,X_{\gamma'})$. The result follows from the first part of the Proposition.
}

\subsection{Characterizations of solutions}

Consider the central extension
\begin{equation}\label{eq:central-extension}
1 \to \Span{J_G} \to G\xrightarrow{\pi} \bar G\to 1
\end{equation}
where $\bar G$ is the quotient group. {For the image $\pi(g)$ of an element $g\in G$ we will simply write $\bar g$.}

{
\Def{\label{def:barNZdG}
Let $\bar N(\ZZ_d,G)$ denote the simplicial set whose $n$-simplices are given by
$$
\bar N(\ZZ_d,G)_n =\set{(\bar g_1,\bar g_2,\cdots,\bar g_n):\, (g_1,\cdots,g_n)\in  N(\ZZ_d,G)_n }.
$$
The simplicial structure maps are similar to that of $N(\ZZ_d,G)$.
}
}

Equivalently, $\bar N(\ZZ_d,G)$ is the orbit space under the action of $N\ZZ_d$ on $N(\ZZ_d,G)$ {which in degree $n$ is given by}
$$
(a_1,\cdots,a_n)\cdot (g_1,\cdots,g_n) = (J^{a_1} g_1,\cdots,J^{a_n} g_n).
$$
This is a free action, hence the   quotient map under this action is a fibration with fiber $N\ZZ_d$:
\begin{equation}\label{eq:fib-N-barN}
N\ZZ_d \to N(\ZZ_d,G) \to \bar N(\ZZ_d,G)
\end{equation}
By the general theory of fibrations it is classified by a cohomology class $\gamma_{G,d}\in H^2(\bar N(\ZZ_d,G))$. More explicitely, this class is the image of the transgression homomorphism in the $E_2$-page of the Serre spectral sequence of the fibration Equation (\ref{eq:fib-N-barN}) (see \cite[Section 6.2]{mccleary2001user}):
$$
d_2:H^1(N\ZZ_d) \to H^2(\bar N(\ZZ_d,G))
$$
Identifying $H^1(N\ZZ_d)\cong \ZZ_d$ we have
\begin{equation}\label{eq:gamma-G-d}
\gamma_{G,d} = d_2(1).
\end{equation}

\begin{lem}\label{lem:triangle-pull-back}
There is a bijection between 
\begin{enumerate}
\item the set of maps
$
f: X\to \bar N(\ZZ_d,G) 
$
such that $f^*(\gamma_{G,d})=[\gamma]$, and
\item the set of commutative diagrams of the form
\begin{equation}\label{diag:triangle}
\begin{tikzcd}
 & N\ZZ_d \arrow[ld,"i"'] \arrow[rd,"\iota"] & \\
X_\gamma \arrow[rr,"\tilde f"] && N(\ZZ_d,G) 
\end{tikzcd}
\end{equation}
\end{enumerate}
\end{lem}
\Proof{
This follows from the classification of fibrations given in Theorem \ref{thm:classification-fib}.
Given $f:X\to \bar N(\ZZ_d,G)$ such that $f^*(\gamma_{G,d})=[\gamma]$ and pulling back $N(\ZZ_d,G)\to \bar N(\ZZ_d,G)$ along $f$ gives a principal $N\ZZ_d$-bundle over $X$ isomorphic to $X_\gamma \to X$. This gives a commutative triangle as in Diagram \ref{diag:triangle}. Conversely, a commutative triangle as in Diagram \ref{diag:triangle} descends to a map $f:X\to \bar N(\ZZ_d,G)$ between the quotient spaces (under the action of $N\ZZ_d$). Comparing the spectral sequences implies that $f^*(\gamma_{G,d})=[\gamma]$.
}

\Cor{ \label{cor:solution-characterization}
Let $(A,b)$ denote the linear system associated to $(X,\gamma)$, where $X$ is a connected simplicial set.
There is a bijective correspondence between the following sets:
\begin{enumerate}
\item the set $\Sol(A,b;G)$ of solutions,
\item the set of group homomorphisms
$
\theta:\Gamma(A,b) \to G
$
such that $\theta(J)=J_G$,
\item the set of group homomorphisms
$
\theta:\pi_1(\ZZ_d,X_\gamma) \to G
$
such that $\theta(e_{1})=J_G$,
\item the set of maps
$
f: X\to \bar N(\ZZ_d,G) 
$
such that $f^*(\gamma_{G,d})=[\gamma]$.
\end{enumerate}
}
\begin{proof}
The bijection between (1) and (3) is proved in Corollary \ref{cor:sol-fund}. {The bijection between (2) and (3) follows from Theorem \ref{thm:pi1-commutative-Xgamma}. 
Part (1) of Proposition \ref{pro:sol-group} gives the bijection between (1) and (2).} Next we prove the bijection between (3) and (4). Proposition \ref{pro:adjunction} implies that there is a bijection between (3) and commutative diagrams given in Diagram \ref{diag:triangle}.
Using Lemma \ref{lem:triangle-pull-back} finishes the proof.
\end{proof}

Known examples of linear systems over $\ZZ_d$ for $d$ an odd integer behave particularly nice.
We state this as a conjecture and provide some evidence for it in Section \ref{sec:linear-system-from-groups}.

\begin{conj}\label{conj:odd}
Let $d>1$ be an odd integer.
Any linear system over $\ZZ_d$ admitting a solution in $G$ also admits a solution in $\ZZ_d$.
\end{conj}

The implications of this conjecture for contextuality is important. 
Let $\chi:G\to U(\hH)$ be a unitary representation.
Proposition \ref{pro:linear-system-contextuality} implies that those linear system admitting a solution in  some group $G$ but not in $\ZZ_d$, produces contextual simplicial distributions for any density operator (see Definition \ref{def:contextuality}).
The conjecture states that those contextual simplicial distributions never arise when $d$ is an odd prime.
However, when $d=2$ this is not the case. For example, the $K_{3,3}$ linear system studied in Example \ref{ex:K33-linear-system} admits a solution in $D_8\ast D_8$ but not in $\ZZ_2$.
This distinction between $d=2$ and $d$ odd has crucial implications in quantum computation \cite{BDBOR17,raussendorf2023role}.

\subsection{Power maps}

For a group $G$ let us write $G_{(d)}$ to denote the set of $d$-torsion elements, i.e., those $g\in G$ satisfying $g^d=1_G$.

\Def{\label{def:Grp-d}
{\rm
Let $\catGrp_{(d)}$ denote the category whose objects are groups and a morphism $G\to H$ is given by a set map $\theta:G_{(d)}\to H$ satisfying
$$
\theta(g_1g_2) = \theta(g_1)\theta(g_2)  
$$
for all $g_1,g_2\in G_{(d)}$ such that $g_1g_2=g_2g_1$.
}
}

Let $\iota:\catGrp \to \catGrp_{(d)}$ denote the inclusion functor.
The nerve space $N(\ZZ_d,\cdot)$ extends to a functor from $\catGrp_{(d)}$ to the category of simplicial sets. 
That is, a morphism $\theta:G\to H$ of the category $\catGrp_{(d)}$ induces a map of simplicial sets $\theta: N(\ZZ_d,G)\to N(\ZZ_d,H)$. This observation follows from Definition \ref{def:commutative-nerve-d}.
There are special morphisms in this category which turn out to be very useful.

\Def{\label{def:power-map}
For an integer $m$, let $\omega_m:G\to G$ denote the set map defined by $g\mapsto g^m$. 
}

In general, the map $\omega_m$ is not a group homomorphism. For instance, $\omega_{-1}$ is a group homomorphism if and only if $G$ is abelian.
However, $\omega_m$ is a morphism of the category $\catGrp_{(d)}$. 
In this case the domain is restricted to the $d$-torsion part $G_{(d)}$.
This map induces a simplicial set map,  denoted by $\omega_m:N(\ZZ_d,G)\to N(\ZZ_d,G)$, which in degree $n$ is given as follows
$$
(g_1,g_2,\cdots, g_n) \mapsto (g_1^m,g_2^m,\cdots, g_n^m).
$$ 
It is straight-forward to verify that this assignment respects the simplicial structure giving us a simplicial set map as claimed.

\Lem{\label{lem:pull-power}
We have $\omega_m^*(\gamma_{G,d})=m\gamma_{G,d}$.
}
\begin{proof}
Raising to the $m$-th power gives a commutative diagram of simplicial sets
\begin{equation}\label{dia:omega-m}
\begin{tikzcd}[column sep=huge,row sep=large]
N\ZZ_d \arrow[r,"{a\, \mapsto ma}"] \arrow[d] & N\ZZ_d \arrow[d]  \\ 
N(\ZZ_d,G) \arrow[r,"\omega_m"] \arrow[d] & N(\ZZ_d,G) \arrow[d] \\
\bar N(\ZZ_d,G) \arrow[r,"\bar\omega_m"] & \bar N(\ZZ_d,G) 
\end{tikzcd}
\end{equation} 
In cohomology the induced map  $H^1(N\ZZ_d)\to H^1(N\ZZ_d)$ is again multiplication by $m$. Recall that the class $\gamma_{G,d}$ is obtained from the transgression $d_2:H^1(N\ZZ_d) \to H^2(\bar N(\ZZ_d,G))$; see Equation (\ref{eq:gamma-G-d}). Comparing the spectral sequences of the fibrations in Diagram (\ref{dia:omega-m}) we obtain the desired result.
\end{proof}

Let us write $d=d_1\cdots d_k$, where $d_i$'s are positive integers such that $\gcd(d_i,d_j)=1$ for all $i\neq j$. 
Let $\gamma_i$ denote the image of $[\gamma]$ under the map $H^2(X,\ZZ_d)\to H^2(X,\ZZ_{d_i})$
induced by the mod $d_i$ reduction homomorphism $\ZZ_d \to \ZZ_{d_i}$.

\begin{lem}\label{lem:mod-di}
Let $q_i=\prod_{j\neq i}d_j$ and $q_i^{-1}$ denote the inverse of $q_i$ in $\ZZ_{d_i}$.
The following diagram of simplicial sets commutes
$$
\begin{tikzcd}[column sep=huge,row sep=large]
N\ZZ_d \arrow[r,"\omega_{q_i}"] \arrow[d] & N\ZZ_{d_i} \arrow[r,"\omega_{q_i^{-1}}"] \arrow[d] & N\ZZ_{d_i} \arrow[d] \\
N(\ZZ_d,G) \arrow[r,"\omega_{q_i}"]  & N(\ZZ_{d_i},G) \arrow[r,"\omega_{q_i^{-1}}"] & N(\ZZ_{d_i},G) 
\end{tikzcd}
$$
Moreover, the composite $\omega_{q_i^{-1}} \circ\omega_{q_i}$ coincides with $N\ZZ_d\to N\ZZ_{d_i}$ induced by the mod $d_i$ reduction.
\end{lem}

\Lem{\label{lem:power-map}
Mod $d_i$ reduction induces an injective map
$$
\set{\,f:X\to \bar N(\ZZ_d,G):\, f^*(\gamma_{G,d})=[\gamma]\,} \longrightarrow \prod_{i=1}^k \set{\,f_i:X\to \bar N(\ZZ_{d_i},G):\, f_i^*(\gamma_{G,d_{i}}) = \gamma_i\,}.
$$  
}
\Proof{ 
Consider the diagram in Lemma \ref{lem:mod-di} and the associated spectral sequence for $H^*(\cdot,\ZZ_{d_i})$. Let $r_i:N(\ZZ_d,G)\to N(\ZZ_{d_i,G})$ denote the composite $\omega_{q_i^{-1}}\circ \omega_{q_i}$ and $\bar r_i$  the induced map $\bar N(\ZZ_d,G)\to \bar N(\ZZ_{d_i},G)$ between the quotient spaces. Combining $\bar r_i$'s gives an injective map
$$
\bar r:N(\ZZ_d,G) \to \prod_i \bar N(\ZZ_{d_i},G).
$$
Given $f:X\to \bar N(\ZZ_d,G)$ we can compose with $\bar r$ and project onto the $i$th factor to obtain $f_i:X\to \bar N(\ZZ_{d_i},G)$. By Lemma \ref{lem:mod-di} the induced map in cohomology
$$
\bar r_i^* : H^2(\bar N(\ZZ_{d_i},G),\ZZ_{d_i}) \to H^2(\bar N(\ZZ_{d},G),\ZZ_{d_i})
$$
sends $\gamma_{G,d_i}$ to the mod $d_i$ reduction of $\gamma_{G,d}$. Therefore $f_i^*(\gamma_{G,d_i})=\gamma_i$ for each $1\leq i\leq k$.
}

In particular, this result can be applied to a prime decomposition of $d$ to reduce problems to the case where $d$ is a prime power. Let $d=\prod_{i=1}^k p_i^{\alpha_i}$ be a prime decomposition. Given a linear system $(A,b)$ over $\ZZ_d$ let us write $(A^{(i)},b^{(i)})$ for the linear system over $\ZZ_{p_i^{\alpha_i}}$ obtained by mod $p_i^{\alpha_i}$ reduction.

\Cor{\label{cor:mod-p-reduction}
Mod $p_i^{\alpha_i}$ reduction reduction induces an injective map between the solution sets:
$$
\Sol(A,b;G) \to \prod_{i=1}^k  \Sol(A^{(i)},b^{(i)};G).
$$
}
\Proof{
Follows from Corollary \ref{cor:solution-characterization} and Lemma \ref{lem:power-map}.
}

Let us turn to the induced map on the algebraic fundamental group of $N(\ZZ_d,G)$ (see Section \ref{sec:CommutativeFund}).   
In general, the power map induces a map on the commutative fundamental group $\pi_1(\ZZ_d,X)$ of an arbitrary simplicial set:
$$
(\omega_m)_*: \pi_1(\ZZ_d,X) \to \pi_1(\ZZ_d,X)
$$
defined by $e_x\mapsto e_x^m$. 
Let $p$ be a positive integer dividing $d$ and  let $q=d/p$.
We can identify $\pi_1(\ZZ_p,X)$ as a subgroup of $\pi_1(\ZZ_d,X)$ via the homomorphism $\iota_q:\pi_1(\ZZ_p,X)\to \pi_1(\ZZ_d,X)$ defined by $e_x\mapsto e_x^q$. 
This homomorphism splits via the composition 
$$
\pi_1(\ZZ_d,X)\xrightarrow{(\omega_q)_*} \pi_1(\ZZ_d,X) \xrightarrow{(\omega_{q^{-1}})_*}
\pi_1(\ZZ_p,X)
$$

Now, we give an application of the $\omega_{-1}$ map to linear systems over $\ZZ_d$ when $d$ is odd. We will write $w$ for a word in a 
finitely presented group, that is, a product of the form $w=e_{x_1} e_{x_2}\cdots e_{x_l}$ where each $e_{x_i}$ is a generator. 
The word $w^\op$ will denote the opposite word given by the product $e_{x_l} \cdots e_{x_2} e_{x_1}$.

\Pro{\label{pro:opposite-odd-Xgamma}
Let $d>1$ be an odd integer.
If the equation $w = e_1^a w^\op$ holds in $\pi_1(\ZZ_d,X_\gamma)$ then $w=w^\op$.
}
\Proof{
The map $\omega_{-1}$ is an automorphism of $\pi_1(\ZZ_d,X_\gamma)$.  
Applying this automorphism to the equation, we obtain that $ w=e_1^{-a} w^\op$, where $ w = e_{x_1}^{-1} e_{x_2}^{-1}\cdots e_{x_l}^{-1}$, also holds. Together with the original equation we obtain
$$
(J^{-a}  w^\op) w =  w (J^a w^\op) \;\Rightarrow\; J^{2a}=1.
$$
Since $d$ is odd, this implies $a=0 \mod d$.
}

Combining this result with Lemma \ref{lem:fund-NZdG} we obtain the following result proved in \cite[Theorem 5]{qassim2020classical} using 
different methods.

\Cor{\label{cor:opposite-odd}
Let $d>1$ be an odd integer.
If the equation $w = J^a w^\op$ holds in $\pi_1 N(\ZZ_d,G)$ then $w=w^\op$.
}
\Proof{
We apply Proposition \ref{pro:opposite-odd-Xgamma} with $X=\bar N(\ZZ_d,G)$ and $[\gamma]=\gamma_{G,d}$. Note that $N(\ZZ_d,G)\cong X_\gamma$.
}

\subsection{The $K_{3,3}$ linear system}\label{sec:K33-linear-system}

In this section we will study the linear system introduced in Example \ref{ex:K33}. 
It will be useful to introduce chains on a simplicial set, a dual notion to cochains.
Given a simplicial set $X$ the set $C_n(X)$  of $n$-chains is the free $\ZZ_d$-module generated by the symbols $[\sigma]$ where $\sigma\in X_n$.
The group $C^n(X)$ of $n$-chains can be identified with the group of $\ZZ_d$-module homomorphisms $C_n(X)\to \ZZ_d$.

Let $(A,b)$ denote this linear system. 
Recall that the dual of the associated simplicial  complex $\Sigma_A$ is given by the complete bipartite graph $K_{3,3}$. 
We can also think of this linear system as one that is obtained from a simplicial set $X$ as depicted in Figure (\ref{fig:K33-b}) and a cocycle $\gamma$ defined by $\sigma\mapsto b_\sigma$.
We will write
$$
[X] = [\sigma_1]+[\sigma_2]+[\sigma_3]-([\sigma_4]+[\sigma_5]+[\sigma_6])
$$
for the $2$-chain. 
We evaluate  $\gamma$ on $[X]$ as follows: $\gamma(X)=\gamma(\sigma_1)+\gamma(\sigma_2)+\gamma(\sigma_3)-(\gamma(\sigma_4)+\gamma(\sigma_5)+\gamma(\sigma_6))$.

\begin{figure}[h!] 
  \centering
  \includegraphics[width=.25\linewidth]{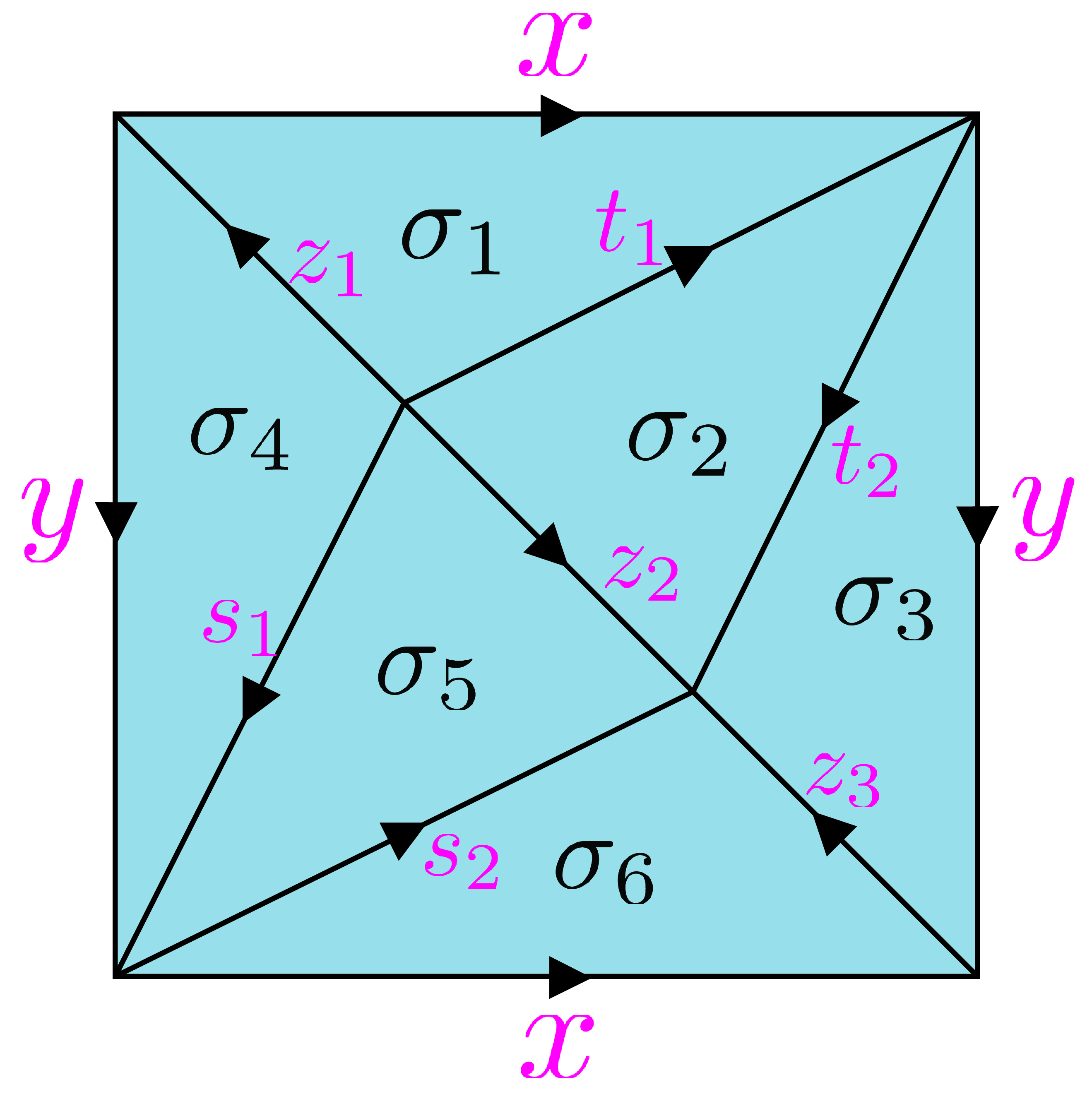}
\caption{
}
\label{fig:K33-fund}
\end{figure} 

\Lem{\label{lem:vanKampen-K33}
For any two nondegenerate  
$1$-simplices $x,y$ that do not belong to a common $2$-simplex the following relation holds in $\pi_1(\ZZ_d,X_\gamma)$:
$$
[e_{x},e_{y}] =  e_1^{\gamma[X]}.
$$ 
}
\Proof{ 
Let $\Aut(K_{3,3})$ denote the automorphism group of $K_{3,3}$. It is well-known that this group is isomorphic to $(\Sigma_3\times \Sigma_3)\rtimes \ZZ_2$  \cite{sreekumar2021automorphism}. 
$\Aut(K_{3,3})$ acts transitively on the set of pairs $\set{x,y}$ that do not belong to a common $\sigma$. This action carries over to the generators in the algebraic fundamental group. Therefore it suffices to consider $x,y$ that live on the boundary of Figure \ref{fig:K33-fund}.
We have
$$
\begin{aligned}
e_x e_y e_x^{-1} e_y^{-1} & =  (e_{z_1}^{-1}e_{t_1}e_1^{\gamma(\sigma_1)})(e_{t_2}e_{z_3}^{-1}e_1^{\gamma(\sigma_3)})(e_{s_2}e_{z_3}^{-1}e_1^{\gamma(\sigma_6)})^{-1}(e_{z_1}^{-1}e_{s_1}e_1^{\gamma(\sigma_4)})^{-1}  \\
&= e_{z_1}^{-1}(e_{t_1}e_{t_2})(e_{s_2}^{-1}e_{s_1}^{-1})e_{z_1} 
e_1^{\gamma(\sigma_1)+\gamma(\sigma_3)-\gamma(\sigma_4)-\gamma(\sigma_6)}\\
&= e_{z_1}^{-1}e_{z_2}e_{z_2}^{-1}e_{z_1} 
e_1^{\gamma(\sigma_1)+\gamma(\sigma_2)+\gamma(\sigma_3)-\gamma(\sigma_4)-\gamma(\sigma_5)-\gamma(\sigma_6)}\\
&= e_1^{\gamma(\sigma_1)+\gamma(\sigma_2)+\gamma(\sigma_3)-\gamma(\sigma_4)-\gamma(\sigma_5)-\gamma(\sigma_6)}\\
&=e_1^{\gamma[X]}.
\end{aligned}
$$ 
}

Recall that under the identification $\Gamma(A,b)\cong \pi_1 (\ZZ_d,X_\gamma)$ we have $e_1\mapsto J$.

\Pro{\label{pro:K33-fund-odd}
When $d>1$ is an odd integer $\pi_1(\ZZ_d,X_\gamma)$ is abelian.  
} 
\Proof{
By Lemma \ref{lem:vanKampen-K33} we have $e_xe_y = J^{\gamma[X]} e_ye_x$. Then Proposition \ref{pro:opposite-odd-Xgamma} implies that $e_xe_y= e_ye_x$.
}
 
\begin{figure}[h!] 
  \centering
  \includegraphics[width=.3\linewidth]{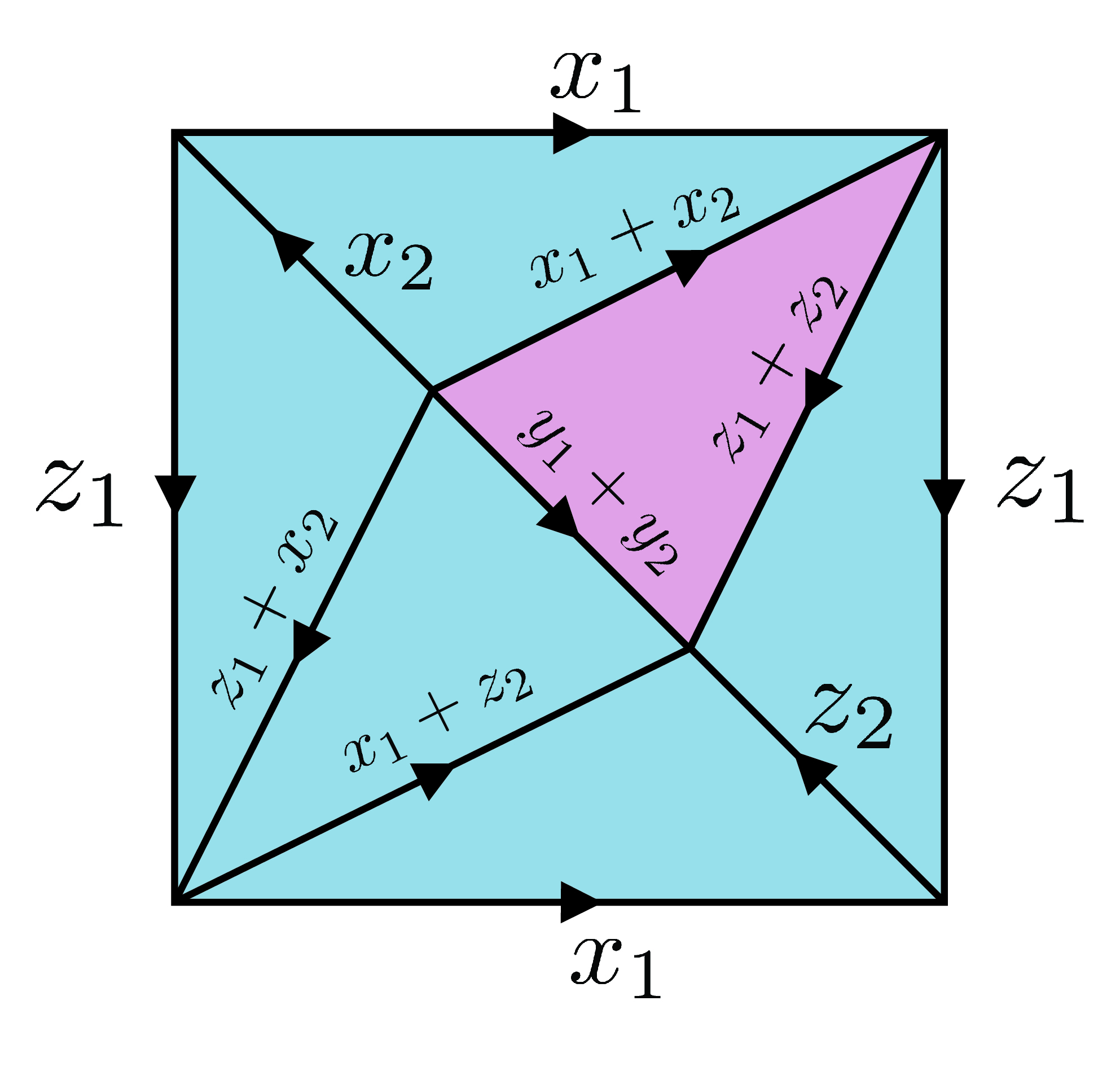}
\caption{
}
\label{fig:D8-D8}
\end{figure}

Let $\gamma$ denote the $2$-cocycle (see Figure (\ref{fig:D8-D8})) defined by
$$
\gamma(\sigma_i) =
\left\lbrace
\begin{array}{ll}
1 & i=2 \\
0 & \text{otherwise.}
\end{array}
\right.
$$ 
The corresponding cohomology class is $[\gamma]=1$.
For the next result we will use the identification $X_\gamma\cong N(\ZZ_2,D_8\ast D_8)$.

\Pro{\label{pro:K33-fund-d=2}
For $d=2$ we have
$$
\pi_1(\ZZ_2,X_\gamma) = \left\lbrace
\begin{array}{ll}
\ZZ_2^{5} & [\gamma]= 0\\
D_8\ast D_8 & [\gamma]= 1.
\end{array}
\right.
$$ 
} 
\begin{proof}
When $[\gamma]=0 $ the result follows from Lemma \ref{lem:vanKampen-K33}. Let $[\gamma]=1$. We can assume that the cocycle representing this class is as given in Figure (\ref{fig:D8-D8}). We rely on the fact that the isomorphism type of the commutative fundamental group only depends on the cohomology class; see Theorem \ref{thm:pi1-commutative-Xgamma}. 
We have a diagram 
$$
\begin{tikzcd}[column sep=huge,row sep=large]
& \Span{e_1}\arrow[d,hook] \arrow[r,"e_1\mapsto -\one"] & \Span{-\one} \arrow[d,hook]\\
K \arrow[r,hook]\arrow[d,equal] & \pi_1(\ZZ_2,X_\gamma) \arrow[d,two heads] \arrow[r,two heads] & D_8\ast D_8 \arrow[d,two heads]\\
K \arrow[r,hook]  &  \pi_1(\ZZ_2,X) \arrow[r,two heads] & \ZZ_2^4  
\end{tikzcd}
$$
where $K$ denotes the kernel of the homomorphism $\pi_1(\ZZ_2,X_\gamma)\to D_8\ast D_8$.
Lemma \ref{lem:vanKampen-K33} implies that $\pi_1(\ZZ_2,X)$ is abelian. From the presentation we see that this group is generated by $e_{x_1}, e_{x_2}, e_{z_1}, e_{z_2}$. On the other hand, it surjects onto $\ZZ_2^4$. Therefore $K=1$. 
\end{proof}

{According to Proposition \ref{pro:linear-system-contextuality} simplicial distributions in the case of $[\gamma]=1$ will be contextual.}
The polytope of simplicial distributions for   $d=2$ is described in \cite{okay2022mermin}.

\section{Linear systems from groups}
\label{sec:linear-system-from-groups}
 
In this section we will study linear systems obtained from groups. 
We begin with a central extension as in (\ref{eq:central-extension}). That is, $G$ is a group with a central element $J_G$ of order $d$ and $\bar G$ is the quotient group by $\Span{J_G}$.
Let 
\begin{equation}\label{eq:phi}
\phi:\bar G \to G
\end{equation}
be a set-theoretic section of the quotient homomorphism $\pi:G\to \bar G$ that preserves the identity element, i.e., $\phi(1_{\bar G})=1_G$.
Given such a section we can construct a $2$-cocycle $\gamma_\phi:G^2\to \ZZ_d$  by 
$$
\gamma_\phi (g,h) = \phi(g)\phi(h)\phi(gh)^{-1}.
$$
It is well-known that the class of the central extension is given by the cohomology class $[\gamma_\phi]\in H^2(\bar G)$; see \cite[Chapter 6]{weibel1995introduction}.
Consider the map of fibrations
\begin{equation}\label{diag:comparison-barN-NbarG}
\begin{tikzcd}[column sep=huge,row sep=large]
N\ZZ_d \arrow[d] \arrow[r,equal] & N\ZZ_d \arrow[d] \\
N(\ZZ_d,G) \arrow[r,hook] \arrow[d] & NG \arrow[d] \\
\bar N(\ZZ_d,G) \arrow[r,hook,"i"] & N\bar G
\end{tikzcd}
\end{equation}
The pull-back of the cochain $\gamma_\phi$, which we denote by $\gamma_{\phi,d}$, represents the cohomology class $\gamma_{G,d}$. {Explicitly, we have $\gamma_{\phi,d}=i_2\circ \gamma_\phi$.}

\Def{\label{def:linear-group}
{\rm
The linear system over $\ZZ_d$ associated to $(G,J_G)$ is the linear system $(A_X,b_\gamma)$ where 
$$
X=\bar N(\ZZ_d,G)\; \text{ and }\; \gamma = \gamma_{\phi,d}.
$$  
We denote this linear system simply by $(A_G,b_\phi)$.
}
}

{
The explicit cocycle $\gamma_{\phi,d}$ can be used to give a twisted product description of $N(\ZZ_d,G)$.

\Pro{\label{pro:NZdG-as-twisted-product}
Let $X=\bar N(\ZZ_d,G)$ and  $\gamma = \gamma_{\phi,d}$.
There is an isomorphism of simplicial sets
$$
X_\gamma\to N(\ZZ_d,G)
$$
} 
\Proof{
In degree $n$ this map is given by 
$$
(a_1,\cdots,a_n) \times (\bar g_1,\cdots,\bar g_n) \mapsto  (J^{a_1}\phi(\bar g_1),\cdots, J^{a_n} \phi(\bar g_n)).
$$
The simplicial structure of the twisted product is such that this assignment gives a well defined simplicial set map.
}

}

\subsection{Homotopical methods}
 
In this section we assume that $G$ is generated by its $d$-torsion elements, i.e., $G=\Span{G_{(d)}}$.
{For details on the homotopical constructions in this section we refer to \cite{O14}.}
 
\Def{\label{def:EZdG} 
Let $E(\ZZ_d,G)$ denote the simplicial set whose set of $n$-simplices is  $G\times N(\ZZ_d,G)_n$ and  simplicial structure maps are given by
$$
d_i(g_0,g_1,\cdots,g_n) = \left\lbrace
\begin{array}{cc}
(g_0,g_1,\cdots,g_ig_{i+1},\cdots, g_n) & 0\leq i <n \\
(g_0,g_1,\cdots,g_{n-1}) & i=n
\end{array}
\right.
$$
and
$$
s_j(g_0,g_1,\cdots,g_n) = s_j(g_0,g_1,\cdots,g_j,1,g_{j-1},\cdots,g_n).
$$
}

There is a simplicial set map $p_d: E(\ZZ_d,G)\to N(\ZZ_d,G)$  
defined in 
degree $n$ by projecting onto the last $n$ coordinates. 
Note that $E(\ZZ_d,G)$ is a simplicial subset of $EG$, total space of the universal principal $G$-bundle $EG\to BG$, whose $n$-simplices are given by $G^{n+1}$ together with a similar simplicial structure.
There is a pull-back diagram of principal $G$-bundles
$$
\begin{tikzcd}[column sep=huge,row sep=large]
E(\ZZ_d,G) \arrow[d,"p_d"] \arrow[r,hook]  & EG \arrow[d,"p"] \\
N(\ZZ_d,G) \arrow[r,hook] & NG
\end{tikzcd}
$$
Note that
\begin{equation}\label{eq:fib-seq-E}
E(\ZZ_d,G) \to N(\ZZ_d,G) \to NG
\end{equation}
is a fibrations sequence. Then the two horizontal maps in the lower part of Diagram (\ref{diag:comparison-barN-NbarG}) can be extended to fibrations:
\begin{equation}\label{diag:extended-barN-NbarG}
\begin{tikzcd}[column sep=huge,row sep=large]
 & N\ZZ_d \arrow[d] \arrow[r,equal] & N\ZZ_d \arrow[d] \\
E(\ZZ_d,G) \arrow[r] \arrow[d,equal] & N(\ZZ_d,G) \arrow[r] \arrow[d] & NG \arrow[d] \\
E(\ZZ_d,G) \arrow[r] & \bar N(\ZZ_d,G) \arrow[r] & N\bar G
\end{tikzcd}
\end{equation}

{Let $|\cdot|$ denote the geometric realization functor.}
Applying $\pi_1(|\cdot|)$ to Diagram (\ref{diag:extended-barN-NbarG}) we obtain a commutative diagram of groups
\begin{equation}\label{diag:groups}
\begin{tikzcd}[column sep=huge,row sep=large]
 & \Span{J} \arrow[r,"{J\mapsto J_G}"] \arrow[d,hook] & \Span{J_G} \arrow[d,hook] &  \\
K(\ZZ_d,G) \arrow[d,equal] \arrow[r,hook] & \Gamma(\ZZ_d,G) \arrow[d,two heads] \arrow[r, two heads] & G \arrow[d,two heads]  \\
K(\ZZ_d,G) \arrow[r,hook] & \bar \Gamma(\ZZ_d,G) \arrow[r,two heads] & \bar G
\end{tikzcd} 
\end{equation} 
where 
$$
K(d,G) = \pi_1|E(\ZZ_d,G)|,\;\; \Gamma(\ZZ_d,G)=\pi_1 N(\ZZ_d,G)\; \text{ and }\; \bar \Gamma(\ZZ_d,G) = \bar N(\ZZ_d,G).
$$ 
 
The group $K(d,G)$ 
can be computed by choosing a maximal tree in the $1$-skeleton; see \cite[Lemma 4]{antolin2021higher} for a similar computation for $\pi_1|E(\ZZ,G)|$. 
The $1$-skeleton is given by a graph whose vertices are the elements of $G$ and edges consist of arrows $g\xrightarrow{g^{-1}h} h$ where $(g^{-1}h)^d=1$.
We choose a maximal tree $T$ as follows: For each vertex $g$ let us choose a sequence of $d$-torsion elements $g_1,\cdots,g_{n(g)}$ satisfying
$
g = g_1\cdots g_{n(g)}.
$
We consider the path from the base point $1\in G$ to the vertex $g$ given by the sequence of arrows  
\begin{equation}\label{eq:path}
p_g:1 \xrightarrow{g_1} g_1 \xrightarrow{g_2} g_1g_2 \xrightarrow{g_3} \cdots \xrightarrow{g_{n(g)}} g_1\cdots g_{n(g)}=g.
\end{equation}
If $g$ is $p$-torsion then we take $n(g)=1$ and $p_g:1\xrightarrow{g} g$.

\Pro{\label{pro:pi1-EZdG}
Let $G$ be a group generated by $d$-torsion elements.
The  
group   
$K(\ZZ_d,G)$
is generated by $e_{g,h}$, where $(g,h)\in G^2$ such that $g^{-1}h\in G_{(d)}$, subject to 
\begin{enumerate}
\item $e_{g,1}=e_{1,g}=1$ for $g\in G_{(d)}$,
\item for $(g,h,k)\in G^3$ such that $(g^{-1}h,h^{-1}k)\in N(\ZZ_d,G)_2$, the relation
$$
e_{g,h} e_{h,k} e_{g,k}^{-1}=1.
$$
\end{enumerate} 
The element $e_{g,h}$ maps to the product $e_{g_1}\cdots e_{g_{n(g)}} e_{g^{-1}h} e_{h_{n(h)}}^{-1}\cdots e_{h_1}^{-1}$ in $\pi_1N(\ZZ_d,G)$.
}
\Proof{ 
The maximal tree $T$ is the union of the paths $p_g$ for each $g\in G$. The generators of the fundamental group are given by the loops based at the identity element. For each $(g,h)\in G$ such that $g^{-1}h \in G_{(d)}$ the generator is given by the loop  
$$e_{g,h}:1\xrightarrow{p_g} g \xrightarrow{g^{-1}h} h \xrightarrow{p_h^{-1}} 1$$ 
The relations come from the degenerate edges and  the $2$-simplices. The former are of the form 
$$
1= e_{s_0(g)} = e_{1,g} \;\;\text{ and }\;\;
1= e_{s_1(g)} = e_{g,1}.
$$
For each triple $(g,h,k)\in G^3$ such that $(g^{-1}h,h^{-1}k)\in N(\ZZ_d,G)_2$, which represents the $2$-simplex $\sigma=(g,g^{-1}h,h^{-1}k)$, we have
$$
1 = e_{d_2 \sigma} e_{d_0 \sigma} e_{d_1\sigma}^{-1} = 
e_{g,h} e_{h,k} e_{g,k}^{-1}.
$$ 
}

\subsection{Finite $p$-groups}

Next we focus on the case where $d$ is given by a prime $p$. 
For a group $H$, we will write $H_\ab$ and $H_\el$ for the largest abelian and elementary abelian $p$-group quotients, respectively.

\Pro{\label{pro:gamma-nonzero-iff-frattini}
Consider the fibration $N\ZZ_p\xrightarrow{i} X_\gamma \to X$.
We have $[\gamma]\neq 0$ if and only if the image of $\pi_1 N\ZZ_p \to \pi_1 X_\gamma$ 
is contained in the kernel of the canonical   homomorphism $\pi_1 X_\gamma \to (\pi_1 X_\gamma)_\el$.
}
\begin{proof}
Consider the five term exact sequence of the Serre spectral sequence of the fibration:
$$
0\to H^1(X) \to H^1(X_\gamma) \xrightarrow{i^*} H^1(N\ZZ_p) \xrightarrow{\delta} H^2(X) \to H^2(X_\gamma)
$$
We have $\delta(1)=[\gamma]\neq 0$. Therefore $i^*$ is the zero map. This is equivalent to the dual map $i_*:H_1(N\ZZ_p) \to H_1(X_\gamma)$ being zero. The result follows from the observation that for a simplicial set $Y$ we have $H_1(Y)=(\pi_1Y)_\el$.
\end{proof}

\Pro{\label{pro:gamma-nonzero-pi1E-zero}
Assume that $[\gamma]\neq 0$ and $K(\ZZ_p,G)=1$. Then $\gamma_{G,p}\neq 0$.
}
\begin{proof}
This follows from Proposition \ref{pro:gamma-nonzero-iff-frattini} and the following diagram
$$
\begin{tikzcd}[column sep=huge,row sep=large]
 & \ZZ_p \arrow[r,equal] \arrow[d,hook] & \ZZ_p \arrow[d,hook]  \\
K(\ZZ_p,G) \arrow[r,hook] & \Gamma(\ZZ_p,G) \arrow[r,two heads] \arrow[d,two heads]  &  G \arrow[d,two heads] \\
& \Gamma(\ZZ_p,G)_\el  \arrow[r,two heads] & G_\el
\end{tikzcd}
$$
Note that only the horizontal sequence of maps is exact.
The composite of the right-hand vertical maps is zero.
 By commutativity of the top square the image of $\ZZ_p\to \Gamma(\ZZ_p,G)$ lands in the 
kernel of $\Gamma(\ZZ_p,G)\to \Gamma(\ZZ_p,G)_\el$.
Hence the class $\gamma_{G,p}$ classifying the fibration $N\ZZ_p \to N(\ZZ_p,G)\to \bar N(\ZZ_p,G)$ is nonzero.
\end{proof}

\Ex{\label{ex:K33-linear-system}
{\rm
{Let $(A,b)$ denote the $K_{3,3}$ linear system studied in Section \ref{sec:K33-linear-system} and $G$ denote the central product $D_8\ast D_8$.}
We will write $X$ for the underlying simplicial set representing the torus depicted in Figure (\ref{fig:K33-fund}). We   identify $\Gamma(A,b)$ with $\pi_1(\ZZ_d,X_\gamma)$ (Theorem \ref{thm:pi1-commutative-Xgamma}).
When $d=2$ Proposition \ref{pro:K33-fund-d=2} implies that if $[\gamma]=1$ then $\pi_1(X_\gamma)\cong G$. 
In this case using Proposition \ref{pro:NZdG-as-twisted-product} we can identify $X_\gamma$ with $N(\ZZ_2,G)$. 
Then we conclude that 
\begin{equation}\label{eq:Gamma-K33}
\Gamma(\ZZ_2,G)\cong G\;\; \text{ and }\;\;K(\ZZ_2,G)=1.
\end{equation}
By Proposition \ref{pro:gamma-nonzero-iff-frattini} we have $\gamma_{G,2}\neq 0$ and by Corollary \ref{cor:solution-characterization} the $K_{3,3}$ linear system does not admit a solution in $\ZZ_2 $. On the other hand,  Equation (\ref{eq:Gamma-K33}) implies that the linear system admits a solution in $G$. For the remaining cases, (1) $d=2$ and $[\gamma]=0$, and (2) $d$ odd, $\Gamma(\ZZ_2,G)$ is abelian, and by Corollary \ref{cor:abelian-sol-group} the linear system admits a solution in $\ZZ_d$.
}
}

\Lem{\cite{AC} \label{lem:odd-p-property}
Let $p$ be an odd prime and $G$ be a nonabelian $p$-group generated by $p$-torsion elements. Then there exists $g,h\in G_{(p)}$ such that $[g,h]\neq 1$ and $g^{-1}h\in G_{(p)}$.
} 
\Proof{
Proof is given in Section \ref{sec:odd-p-property}.
}

\Cor{\label{cor:odd-p-K}
Let $p$ be an odd prime and $G$ be a nonabelian $p$-group generated by $p$-torsion elements. Then $K(\ZZ_p,G)$ is nontrivial.
} 
\Proof{
Let $g,h$ be $p$-torsion elements satisfying $[g,h]\neq 1$ and $g^{-1}h$ is $p$-torsion, obtained by  Lemma \ref{lem:odd-p-property}.
Then Proposition \ref{pro:pi1-EZdG} implies that $e_{g,h}$ is a non-trivial element in $K(\ZZ_p,G)$.
}

Assume that the element $e_{g,h}$ constructed in Corollary \ref{cor:odd-p-K} is central in $\Gamma(\ZZ_p,G)$.
Consider the diagram of group extensions
\begin{equation}\label{dia:implies-conj}
\begin{tikzcd}[column sep=huge,row sep=large]
\Span{e_{g,h}} \arrow[d,hook] \arrow[r,hook] & K(\ZZ_p,G) \arrow[d,hook] \\
\tilde G \arrow[d,two heads] \arrow[r,hook] &  \bar \Gamma(\ZZ_p,G) \arrow[d,two heads] \\
\bar G \arrow[r,equal] & \bar G
\end{tikzcd}
\end{equation} 
If the generator $1\in H^1(\Span{e_{g,h}})\cong \ZZ_p$ hits the extension class $\gamma_G$ then
$\gamma_{G,p}=0$.
As we will see in the next section this is indeed the case for extraspecial $p$-groups when $p>2$.
{Then using Corollary \ref{cor:solution-characterization}, we see that 
this observation provides some evidence for  Conjecture \ref{conj:odd} when $d$ is an odd prime. 
}

\subsection{Extraspecial $p$-groups}
\label{sec:extraspecial}

A $p$-group $E$ is called extraspecial if $Z(E)=[E,E]=\Phi(E)\cong \ZZ_p${; see \cite{As86}}.
The order of an extraspecial $p$-group is given by $p^{2n+1}$ where $n\geq 1$.
When $n=1$ there are two types $E_1^+$ and $E_1^-$. 
For $p=2$, $E_1^+$ is the dihedral group $D_8$ and $E_1^+$ is the quaternion group $Q_8$. 
For $p>2$, $E_1^+$ is the group of triangular $3\times 3$ matrices over $\ZZ_p$ with $1$'s on the diagonal, and $E_1^-$ is the semidirect product of a cyclic group of order $p^2$ by a cyclic group of order $p$ acting non-trivially. 
In general, for $n\geq 1$ we have  
$$
E^+_n = \underbrace{E_1^+\ast \cdots \ast E_1^+}_n \text{ and } E^-_n = \underbrace{E_1^+ \ast \cdots \ast E_1^+}_{n-1}\ast E_1^-
$$
where $G\ast H$ denotes the central product.
Let $J_E$ denote a  generator of $Z(E)$. 
An extraspecial $p$-group fits in a central extension
$$
1\to \Span{J_E} \to E_n \to \ZZ_p^{2n} \to 0. 
$$ 
Finally, for $p=2$ we will also consider almost extraspecial $2$-groups, which are defined by $E_n^0 = E_n^+ \ast \ZZ_4$. In this case the central product is with respect to the cyclic subgroup of order $2$ in $\ZZ_4$. 
 
When $p>2$, the group $E_n^-$ is not generated by $p$-torsion elements; hence will be omitted from our consideration. 
In the next result we compute $\pi_1 N(\ZZ_p,E)$ for the remaining types of extraspecial $p$-groups. 
We will need the following construction:
For a group $G$ let us write $\hat G$ for the subgroup of $G\times G$ generated by $(g,g^{-1})$.

\Thm{\label{thm:extraspecial}
For an extraspecial $p$-group $E$ the group $\Gamma(\ZZ_,E)$ is described as follows:
\begin{enumerate}
\item For $p=2$,  
$\Gamma(\ZZ_2,E)\cong E$ if $E=E_n^+$ and $n\geq 2$, or $E=E_n^-$ and $n\geq 3$.

\item For an almost extraspecial $2$-group, $\Gamma(\ZZ_2,E_n^0)\cong \hat E_n^0$ if $n\geq 2$.

\item For $p>2$, $\Gamma(\ZZ_2,E_n^+)\cong \hat E_n^+$.
\end{enumerate}
}
\begin{proof}
Let us begin with the case $p=2$. 
Observe that the linear system associated to $(E_2^+,J_E)$ can be identified with the odd parity ($[\gamma]=1$) linear system associated to  $K_{3,3}$. More precisely, $X_\gamma \cong N(\ZZ_2,E_2^+)$. Therefore Proposition \ref{pro:K33-fund-d=2} implies that the statement holds for $E_2^+$.
To generalize this to $n\geq 2$, the argument is similar. Consider the diagram of groups 
$$
\begin{tikzcd}[column sep=huge,row sep=large]
& \Span{J}\arrow[d,hook] \arrow[r,"J\mapsto J_E"] & \Span{J_E} \arrow[d,hook]\\
K(\ZZ_2,E) \arrow[r,hook]\arrow[d,equal] & \Gamma(\ZZ_2,E) \arrow[d,two heads] \arrow[r,two heads] & E \arrow[d,two heads]\\
K(\ZZ_2,E) \arrow[r,hook]  &  \bar\Gamma(\ZZ_2,E) \arrow[r,two heads] & \ZZ_2^{2n}  
\end{tikzcd}
$$
The key fact we need is the following: Given two pairs of $p$-torsion elements $(g_1,g_2)$ and $(h_1,h_2)$ of $E$ satisfying $[g_1,g_2]=[h_1,h_2]$ there exists an automorphism $\phi$ of $E$ fixing the central element $J_E$ such that $\phi(g_1)=h_1$ and $\phi(g_2)=h_2$.
This is a well-known fact and follows from Witt's lemma \cite[Chapter 7]{As86}. Using the action of the automorphism group of $E$ on $\Gamma(\ZZ_2,E)$ and Lemma \ref{lem:vanKampen-K33} we see that any pair of $2$-torsion elements satisfying $[g_1,g_2]=J_E$ will satisfy $[e_{g_1},e_{g_2}]=J$. 
Then $\bar\Gamma(\ZZ_2,E)$ turns out to be isomorphic to $\ZZ_2^{2n}$. 
This implies that $K(\ZZ_2,E)=1$. 
A similar argument works for $E=E_n^-$ when $n\geq 3$. 

For the almost extraspecial group, the simplicial set $N(\ZZ_2,E_n^0)$ can be identified with $N(\ZZ_2,E_n^+)$. This follows from the basic fact that there is a bijection between the partially ordered set (poset) of the elementary abelian subgroups of $E_n^0$ and the poset of abelian subgroups of $E_n^+$; see \cite[Chapter 8]{As86}. 
The fundamental group of the latter simplicial set is computed in \cite{O15}.

For $p>2$ note that $N(\ZZ_p,E)$ coincides with $N(\ZZ,E)$ since every element of $E$ is $p$-torsion. 
The fundamental group of $N(\ZZ,E)$ for $n\geq 2$  is isomorphic to $\hat E$ as shown in \cite{O15}.
\end{proof}

\Cor{\label{cor:extraspecial-p-odd}
If a linear system admits a solution in $E_n^+$ for $p>2$ then it admits a solution in $\ZZ_p$.
}
\Proof{
By Theorem \ref{thm:extraspecial} we have
$$
\bar \Gamma(\ZZ_p,E) \cong \hat E / \Span{(J_E,J_E^{-1})}.
$$
This group turns out to be isomorphic to $E_n^+$. Therefore in the Diagram (\ref{dia:implies-conj}) the extension class $\gamma_E$ is hit by the transgression map. Thus $\gamma_{E,p}=0$.
}

\subsection{Higher odd prime torsion groups}

In this section we give more examples in support of Conjecture \ref{conj:odd}. For details see \cite{frembs2022no}. 
Let $\hH=\CC\ZZ_p$ denote the vector space with basis $\set{b_a:\,a\in \ZZ_p}$. 
For $m\geq 1$, we define a subgroup of the special unitary group $\SU(\hH)$:
$$
E_1(p^m) = \Span{T_{(p^m)},X}
$$
where $T$ is the maximal torus given by the diagonal matrices in $\SU(\hH)$ and $X$ is the permutation matrix defined by $X b_a =b_{a+1}$.
We define
$$
E_n(p^m) = \underbrace{E_1(p^m) \otimes \cdots \otimes E_1(p^m)}_n.
$$
Note that for $m=1$ we obtain $E_n^+$. 

An element of $T$ can be described by a function $\xi:\ZZ_p\to U(1)$ such that $\prod_{q\in \ZZ_p} \xi(q)=1$. The corresponding diagonal matrix will be denoted by $D(\xi)$.
The function $\xi$ can be expressed as
$$
\xi(q) = e^{\sum_{j=1}^m  2\pi i f_j(q)/p^j } 
$$
where $f_j(q) = \sum_{b\in \ZZ_p} \nu_{j,b} q^b$ and $\nu_{j,b}\in \ZZ_p$.
We define a function
\begin{equation}\label{eq:phi}
\phi: E_n(p^m) \to E_n(p)
\end{equation}
as follows
$$
\phi(M_1\otimes \cdots \otimes M_n) = \phi_1(M_1)\otimes \cdots \otimes \phi_1(M_n)
$$
where for $M=D(\xi)X^b$ we have
$$
\phi_1(M) = \left\lbrace
\begin{array}{ll}
D(\omega^{\nu_{1,0}+\nu_{2,0}+\nu_{1,1}q}) & b=0 \\
D(\omega^{\nu_{1,0}+\nu_{1,1}q})X & b=1 \\
\phi_1(M^{b^{-1}})^{b} & 1<b<p.
\end{array}
\right.
$$

\Thm{ \label{thm:Markus-paper}
Let $p$ be an odd prime.
The function $\phi$ defined in Equation (\ref{eq:phi}) is a morphism in $\catGrp_{(p)}$ that splits the inclusion $E_n(p)\to E_n(p^m)$. As a consequence the induced homomorphism
$$
\Gamma(\ZZ_p,E_n(p)) \to \Gamma(\ZZ_p,E_n(p^m)) 
$$
splits in $\catGrp$.
}
\Proof{
Proof is given in \cite[Theorem 2]{frembs2022no}. 
}

\Cor{\label{cor:Markus-solution}
If a linear system admits a solution in $E_n(p^m)$ then it admits a solution in $\ZZ_p$.
}
\Proof{
Follows from Theorem \ref{thm:Markus-paper}, Corollary \ref{cor:extraspecial-p-odd} and Corollary \ref{cor:abelian-sol-group}.
}

\bibliography{bib.bib}
\bibliographystyle{alpha}

\appendix

\section{Proof of Proposition \ref{pro:reduction}} \label{sec:proof-reduction}

\Lem{\label{lem:delta-b}
Under the connecting homomorphism $\delta$ in (\ref{eq:connecting-hom}) the class $b$ maps to $\gamma_b$.
}
\Proof{
We will use the definition of the connecting homomorphism to construct a $2$-cocycle representing $\delta(b)$. This is done in two steps \cite[Chapter 6]{weibel1995introduction}: (1) First $b$ is lifted to a $1$-cochain on $N(\ZZ_d,\Sigma)$. (2) Then the coboundary $d$ of $N(\ZZ_d,\Sigma)$ is applied.
We have the following identification
$$
\begin{aligned}
N(\ZZ_d,\Sigma)_1 &=  \set{s:\Sigma_0\to \ZZ_d:\, \supp(s)\in \Sigma} \\
&\cong  \frac{\set{s:\Sigma_0\to \ZZ_d:\, \supp(s)\in \Sigma}}{aA_i \sim 0} \vee \bigvee_{i=1}^r \set{aA_i:\, a\in \ZZ_d} \\
&\cong \bar N(\ZZ_d,\Sigma)_1 \vee \bigvee^r (N\ZZ_d)_1.
\end{aligned}
$$
The group of $1$-cochains $C^1(\vee_i N\ZZ_d)$ is given by $\ZZ_d^{(\vee_i N\ZZ_d)_1}= \ZZ_d^{\sqcup_i \ZZ_d} \cong \prod_i \ZZ_d^{\ZZ_d}$ and the isomorphism in Equation (\ref{eq:H1-N}) can be obtained by identifying $H^1(\vee_i N\ZZ_d)$ with the subset $\prod_i \catGrp(\ZZ_d,\ZZ_d)$ of group homomorphisms. Consider the exact sequence of $1$-cochains
$$
0 \to \ZZ_d^{\bar N(\ZZ_d,\Sigma)_1} \to \ZZ_d^{ N(\ZZ_d,\Sigma)_1} \to \ZZ_d^{\sqcup_i \ZZ_d} \to 0.
$$
The lifting $\tilde b$ we choose is the one given in Equation (\ref{eq:tilde-b}), which is defined to be zero for all elements in $N(\ZZ_d,\Sigma)_1-(\vee_i N\ZZ_d)_1$. Then $\gamma_b$ is obtained by applying the coboundary $d$.
}

\begin{proof}[{\bf Proof of Proposition \ref{pro:reduction}}]
Given two functions $r,s:\Sigma_0\to \NN$ we will write $r\leq s$ if $r(v)\leq s(v)$ for all $v\in \Sigma_0$. This induces a partial order on the set $\NN^{\Sigma_0}$. Throughout we think of a function $s:\Sigma_0\to \ZZ_d$ as an element of $\NN^{\Sigma_0}$ by identifying $\ZZ_d$ with the subset $\set{0,1,\cdots,d-1}\subset \NN$. 
Then $N(\ZZ_d,\Sigma)_1$ is regarded as a subset of $\NN^{\Sigma_0}$ and $\leq$  induces a partial order on this subset.
Recall that we have
$$
N(\ZZ_d,\Sigma)_1 = \bar N(\ZZ_d,\Sigma) \vee A
$$
where $A = \vee_{i=1}^{r} \set{aA_i:\, a \in \ZZ_d}$.  
The group  $\Gamma(A_X,b_\gamma)$ is generated by $e_s$ where $s\in \bar N(\ZZ_d,\Sigma)$ and the product relations come from $2$-simplices of $\bar N(\ZZ_d,\Sigma)$: For $(s_1,s_2) \in \bar N(\ZZ_d,\Sigma)_2$ the corresponding relation is
\begin{equation}\label{eq:rel}
e_{s_1} e_{s_2} e_{s_1+s_2}^{-1}  = J^{-\gamma(s_1,s_2)}. 
\end{equation}
See Definition \ref{def:solution-group} and \ref{def:linear-system-of-simplicial-set}.
For $v\in \Sigma_0$ let us write $\delta^v:\Sigma_0\to \ZZ_d$ for the delta function
\[
\delta^v(v')=\left\{
\begin{array}{ll}
1&\textrm{if }v'=v\\
0&\textrm{otherwise.}
\end{array}
\right.
\]
First we observe that $\Gamma(A_X,b_\gamma)$ is generated by $J$ and $e_{\delta^v}$ where $v\in \Sigma_0$, i.e., every generator of the form $e_s$ can be decomposed as a product of those.
This can be proved by   induction over the partial order. The claim holds for the minimal non-trivial functions, that is, the delta functions $\delta^v$. 
For an arbitrary $s\notin A$ we can write $s= r+r'$ for some $r,r'<s$.
Then using Eq.~(\ref{eq:rel}) we have
\begin{equation}\label{eq:decomposition}
e_s = 
\left\lbrace
\begin{array}{ll}
 J^{\alpha}  & r,r'\in A\\
 e_{r'} J^{\alpha}  & r\in A,\; r' \notin A\\
 e_r e_{r'} J^{\alpha}  & r,r'\notin A
\end{array}
\right.
\end{equation}
where $\alpha$ can be computed using the formula for $\gamma(r,r')$. By induction $e_r$ and $e_{r'}$ is a product of $J$ and $e_{\delta^v}$'s.  

We will show that the following maps:
$$ 
\phi: \Gamma(A,b) \to \Gamma(A_X,b_\gamma) 
$$
defined by $\phi(e_v)=e_{\delta^v}$ and
$$
\psi: \Gamma(A_X,b_\gamma)  \to  \Gamma(A,b)
$$
defined by   $\psi(e_{\delta^v})=e_v$ and $J\mapsto J$ in both cases, give well-defined group homomorphisms. These maps clearly respect the $d$-torsion and commutativity relations. It suffices to show that they also respect the product relations. This will prove the isomorphism as these maps are inverse of each other.

We will say that $s$ is {\it reduced} if there exists no $r\in A$ such that $r\leq s$. We will use the following relation
\begin{equation}\label{eq:union}
\prod_{v\in \supp(s_1+s_2)} e_{\delta^v}^{(s_1+s_2)(v)} = \prod_{v\in \supp(s_1)\,\cup\,(s_2)} e_{\delta^v}^{(s_1+s_2)(v)} =  
\prod_{v\in \supp(s_1)} e_{\delta^v}^{s_1(v)}
\prod_{v\in \supp(s_2)} e_{\delta^v}^{s_2(v)}.
\end{equation}
Note that we regard $s_1,s_2$ as functions with target $\NN$.
First we will show that $\phi$ respects the product relations.
To this end, given $aA_i\neq 0$ we will show that the following relation holds in $\Gamma(A_X,b_\gamma)$:
\begin{equation}\label{eq:aAi}
\prod_{v\in {\supp(aA_i)}} e_{\delta^v}^{a A_{i}(v)} = J^{ab_i}.
\end{equation}

\begin{itemize}
\item Let $aA_i\neq 0$ be a minimal element of $A$.
Choose $v\in \supp(aA_i)$ and define $s=aA_i-\delta^v$.
The function $s$ is reduced, thus Eq.~(\ref{eq:decomposition}) gives
$$
\begin{aligned}
J^{ab_i} & = e_s e_{\delta^v}\\
& =  \left( \prod_{ w\in \supp(s)} e^{s(w)}_{\delta^w} \right)  e_{\delta^v} \\
& =  \prod_{w\in \supp(aA_i)} e^{aA_i(w)}_{\delta^w} . 
\end{aligned}
$$
\item For an arbitrary $a A_i$, let $bA_j$ be a predecessor in $A$, i.e, if $b A_j \leq r \leq aA_i$ then either $r=a A_i$ or $r=bA_j$. 
Let $s= aA_i-bA_j$.
By induction and using Eq.~(\ref{eq:decomposition}) and Eq.~(\ref{eq:union}) we have
$$
\begin{aligned}
J^{ab_i} & = J^{ab_i-bb_j} J^{bb_j} \\
&= e_s \left( \prod_{v\in \supp(bA_j)} e_{\delta^v}^{bA_j(v)} \right) \\
&= \left( \prod_{w\in \supp(s)} e_{\delta^w}^{s(w)} \right)  \left( \prod_{v\in \supp(bA_j)} e_{\delta^v}^{bA_j(v)} \right) \\
&=   \prod_{w\in \supp(aA_i)} e_{\delta^w}^{aA_i(w)}. 
\end{aligned}
$$ 
\end{itemize}
This proves that $\phi$ is a group homomorphism.

Let us turn to $\psi$. 
We will need the following counterpart of Eq.~(\ref{eq:union}):
\begin{equation}\label{eq:union-supp-plus}
\prod_{v\in \supp(s_1+s_2)} e_{v}^{(s_1+s_2)(v)} =  
\prod_{v\in \supp(s_1)} e_{v}^{s_1(v)}
\prod_{v\in \supp(s_2)} e_{v}^{s_2(v)}.
\end{equation}
We begin by observing that 
\begin{equation}\label{eq:psi-e-s}
\psi(e_s) = \prod_{v\in \supp(s)} e_v^{s(v)}
\end{equation}
for an arbitrary function $s\notin A$. This can be proved using Eq.~(\ref{eq:decomposition}) and induction on the partial order. It is clear that this holds for a reduced $s$. Otherwise, $s=r+r'$ where $r\in A$ and $r'<s$. Then by induction and using Eq.~(\ref{eq:decomposition}) and Eq.~(\ref{eq:union-supp-plus}) we have:
\begin{itemize}
\item If $r'\in A$ then
$$
\begin{aligned}
\psi(e_s) &= \psi(J^{\gamma(r,r')})\\
&= J^{ab_i + bb_j} \\
& = \prod_{v\in \supp(r)} e_{v}^{r(v)} 
\prod_{v\in \supp(r')} e_{v}^{r'(v)}\\
& = \prod_{v\in \supp(r+r')} e_{v}^{r(v)+r'(v)} \\   
&= \prod_{v\in \supp(s)} e_v^{s(v)}
\end{aligned}
$$
where $r'=bA_j$.
\item $r'\notin A$
$$
\begin{aligned}
\psi(e_s) &= \psi(e_{r'}J^{\gamma(r,r')})\\
&= \psi(e_{r'})  J^{ab_i} \\
&= \prod_{v\in \supp(r')} e_{v}^{r'(v)} \prod_{v\in \supp(r)} e_{v}^{r(v)} \\
&= \prod_{v\in \supp(s)} e_v^{s(v)}.
\end{aligned}
$$
\end{itemize}
To finish the proof we will use Eq.~(\ref{eq:union-supp-plus}) and Eq.~(\ref{eq:psi-e-s}).
There are eight cases:
\begin{itemize}
\item $s_1,s_2,s_1+s_2\notin A$: We have
$$
\begin{aligned}
\psi(e_{s_1} e_{s_2} e_{s_1+s_2}^{-1}) &= \prod_{v\in \supp(s_1)} e_v^{s_1(v)} \prod_{v\in \supp(s_2)} e_v^{s_2(v)} \prod_{v\in \supp(s_1+s_2)} e_v^{-(s_1(v)+s_2(v))} \\
&= \prod_{v\in \supp(s_1)} e_v^{s_1(v)} \prod_{v\in \supp(s_2)} e_v^{s_2(v)} \prod_{v\in \supp(s_1)} e_v^{-s_1(v)} \prod_{v\in \supp(s_2)} e_v^{-s_2(v)} \\
&=1.
\end{aligned}
$$ 
Note that we use here the fact that $(s_1,s_2)$ is a $2$-simplex in $\bar{N}(\ZZ_d,G)$, thus $\supp(s_1)\cup\supp(s_2)$ belongs to $\Sigma$.
Thus there is a maximal simplex $\sigma_i\in \Sigma$ that contains $\supp(s_1)$, $\supp(s_2)$, and $\supp(s_1+s_2)$. Thus by commutativity relation (Definition \ref{def:solution-group}) all elements $e_v$ in the equation above commute. A similar observation applies to the remaining cases.

\item $s_1,s_2\notin A$ and $s_1+s_2\in A$: We have
$$
\begin{aligned}
\psi(e_{s_1} e_{s_2} J^{\gamma(s_1,s_2)}) &= \prod_{v\in \supp(s_1)} e_v^{s_1(v)} \prod_{v\in \supp(s_2)} e_v^{s_2(v)} J^{\gamma(s_1,s_2)}\\
    &= \prod_{v\in \supp(s_1+s_2)} e_v^{s_1(v)+s_2(v)} J^{\gamma(s_1,s_2)}\\
       &= \prod_{v\in \supp(aA_i)} e_v^{aA_i(v)} J^{\gamma(s_1,s_2)}\\
&= J^{ab_i} J^{-ab_i} =1       
\end{aligned}
$$
where $s_1+s_2=aA_i$.

\item $s_1\in A$ and $s_2,s_1+s_2\notin A$: We have
$$
\begin{aligned}
\psi( e_{s_2} e_{s_1+s_2}^{-1}J^{\gamma(s_1,s_2)}) &=  \prod_{v\in \supp(s_2)} e_v^{s_2(v)} \prod_{v\in \supp(s_1+s_2)} e_v^{-(s_1(v)+s_2(v))} J^{\gamma(s_1,s_2)} \\
&=  \prod_{v\in \supp(s_2)} e_v^{s_2(v)} \prod_{v\in \supp(s_1)} e_v^{-s_1(v)}\prod_{v\in \supp(s_2)} e_v^{-s_2(v)} J^{\gamma(s_1,s_2)}  \\
&=  J^{-ab_i} J^{ab_i}=1 
\end{aligned}
$$
where $s_1=aA_i$. The case  $s_2\in A$ and $s_1,s_1+s_2\notin A$ is similar.

\item $s_1,s_1+s_2\in A$ and $s_2\notin A$: We have
$$
\begin{aligned}
\psi(  e_{s_2} J^{\gamma(s_1,s_2)}) &= \prod_{v\in \supp(s_2)} e_v^{s_2(v)} J^{\gamma(s_1,s_2)}\\
&= \prod_{v\in \supp(s_1)} e_v^{-s_1(v)} \prod_{v\in \supp(s_1)} e_v^{s_1(v)}  \prod_{v\in \supp(s_2)} e_v^{s_2(v)} J^{\gamma(s_1,s_2)}\\
&= \prod_{v\in \supp(s_1)} e_v^{-s_1(v)} \prod_{v\in \supp(s_1+s_2)} e_v^{s_1(v)+s_2(v)}  J^{\gamma(s_1,s_2)}\\
&= J^{-ab_i} J^{bb_j} J^{ab_i-bb_j}=1 
\end{aligned}
$$
where $s_1=aA_i$ and $s_1+s_2=bA_j$. The case $s_2,s_1+s_2\in A$ and $s_1\notin A$ is similar.

\item $s_1,s_2\in A$ and $s_1+s_2\notin A$: We have
$$
\begin{aligned}
\psi( e_{s_1+s_2}^{-1}J^{\gamma(s_1,s_2)}) &=
\prod_{v\in \supp(s_1+s_2)} e_v^{-(s_1(v)+s_2(v))} J^{\gamma(s_1,s_2)} \\
&=\prod_{v\in \supp(s_1)} e_v^{-s_1(v)}\prod_{v\in \supp(s_2)} e_v^{-s_2(v)} J^{\gamma(s_1,s_2)} \\
&=J^{-aA_i} J^{-bA_j} J^{ab_i+bb_j}=1
\end{aligned}
$$
where $s_1=aA_i$ and $s_2=bA_j$.

\item $s_1,s_2,s_1+s_2 \in A$: We have
$$
\begin{aligned}
\psi(J^{\gamma(s_1,s_2)}) &= J^{\gamma(s_1,s_2)} \\
&=\left( \prod_{v\in \supp(s_1)} e_v^{s_1(v)} \right)^{-1} \prod_{v\in \supp(s_1)} e_v^{s_1(v)}  
\left( \prod_{v\in \supp(s_2)} e_v^{s_2(v)} \right)^{-1} \prod_{v\in \supp(s_2)} e_v^{s_2(v)} J^{\gamma(s_1,s_2)} \\
&=   \left( \prod_{v\in \supp(s_1)} e_v^{s_1(v)} \right)^{-1}\left( \prod_{v\in \supp(s_2)} e_v^{s_2(v)} \right)^{-1} \prod_{v\in \supp(s_1+s_2)} e_v^{s_1(v)+s_2(v)} J^{\gamma(s_1,s_2)} \\
&= J^{-aA_i} J^{-bA_j} J^{cA_k} J^{aA_i-cA_k+bA_j}=1
\end{aligned}
$$
where $s_1=aA_i$, $s_2=bA_j$ and $s_1+s_2 = cA_k$.
\end{itemize}  
This proves that $\psi$ is a group homomorphism. By definition $\phi$ is the inverse. Thus we obtain that the two groups are isomorphic.   
\end{proof}

\section{Classification of fibrations}\label{sec:ClassFib}

Let $X$ be a simplicial set and $[\gamma]$
 be a cohomology class in $H^2(X)$ represented by a normalized cocycle $\gamma:X_2\to \ZZ_d$ {(see Equation (\ref{eq:normalized}))}. 
Using the cocycle we can define a simplicial group homomorphism $f_{\gamma}: GX \to N\ZZ_d$, where $GX$ denotes Kan's loop space.
First observe that a simplicial group homomorphism $f:GX\to N\ZZ_d$ is uniquely specified by a group homomorphism $f_1:F[X_2-s_0X_1]\to \ZZ_d$ that satisfies 
$$
\begin{aligned}
f_1(d_1[\tau]) &= f_1(d_2[\tau]) f_1(d_0[\tau]) &\Leftrightarrow \\
f_1([d_2\tau])  &= f_1([d_3\tau]) f_1([d_1\tau][d_0\tau]^{-1}) &\Leftrightarrow\\
df_1(\tau) &=0
\end{aligned} 
$$
for all $\tau \in X_3-s_0 X_2$.
Here in the second line we used the face maps of $GX$ and in the third line we regard $f_1$ as a function $(X_2-s_0 X_1)\to \ZZ_d$. Therefore given a cocycle $\gamma$ the map $f_{\gamma}$ is defined by
$$
(f_{\gamma})_1([\sigma]) = \gamma(\sigma).
$$
In dimensions $n\geq 2$ the map is defined as follows: Let $e_i:(GX)_n\to (GX)_1$  denote the $i$-th spine operator. For $\tau\in X_{n+1}-s_0X_n$ we have
$$
f_{\gamma}([\tau]) = (f_1(e_1[\tau]), f_1(e_2[\tau]),\cdots, f_1(e_n[\tau]) ).
$$
Using the adjoint pair $G\adjoint \overline W$ we obtain a simplicial set map
$g: X\to \overline W(N\ZZ_d)$. 
In dimension $n$ this map is given by 
$$
g(\tau)=(f([\tau]),d_0 f([\tau]),\cdots ,d_0^{n-1} f([\tau])) .
$$
We define $X_\gamma$ to be the pull-back of the universal bundle
$$
\begin{tikzcd}[column sep=huge,row sep=large]
X_\gamma \arrow[d] \arrow[r] &  W(N\ZZ_d) \arrow[d] \\
X \arrow[r,"g"] & \overline W(N\ZZ_d)
\end{tikzcd}
$$
We can provide a description of $X_\gamma$ as a twisted cartesian product by pulling-back the canonical twisting function of the universal bundle, which is given by projection onto the first coordinate $\overline W(N\ZZ_d)_n \to (N\ZZ_d)_{n-1}$. 
For $n=2$ the twisting function $\eta:X_2\to \ZZ_d$ is given by 
$$
\eta(\sigma) = \gamma(\sigma)
$$
and for $n\geq 3$ the twisting function $\eta:X_n\to \ZZ_d^{n-1}$ can be computed using $g$:
\begin{equation}\label{eq:eta}
\eta(\tau) = f([\tau])= (f_1(e_1[\tau]),\cdots,f_1(e_{n-1}[\tau]) ).
\end{equation} 
The twisted cartesian product $N\ZZ_d\times_\eta X$ is defined as follows:
\begin{itemize}
\item the set of $n$-simplices are given by $\ZZ_d^n\times X_n$,
\item the simplicial structure maps are given by
$$
d_i(\alpha , \tau) = \left\lbrace
\begin{array}{ll}
(\eta(\tau) +d_0(\alpha), d_0\tau) & i=0 \\
(d_i\alpha,d_i\tau) & \text{otherwise,}
\end{array}
\right. 
$$
and $s_j(\alpha,\tau)=(s_j\alpha,s_j\tau)$. 
\end{itemize}
We will write $N\ZZ_d \times_{\gamma} X$ for the twisted product $N\ZZ_d\times_\eta X$ to indicate that it comes from a $2$-cocycle.

\Thm{\label{thm:classification-fib}
There is a bijective correspondence between the set of isomorphism classes of principal $N\ZZ_d$-bundles over $X$ and the cohomology group $H^2(X)$. 
}
\Proof{
This is a particular case of the classical result on the classification of principal $K$-bundles for a simplicial group $K$; see \cite[Chapter IV]{May67}.
The idea is as follows: Given $[\gamma]$ 
one associates the twisted cartesian product $X_\gamma= N\ZZ_d \times_{\gamma} X$ constructed above. Conversely, back-tracking the same argument produces a cohomology class from a twisted cartesian product. Up to isomorphism principal bundles are given by twisted cartesian products and changing a cocycle by a coboundary produces an isomorphic twisted cartesian product. Combining these facts gives the classification theorem. 
}

\section{Proof of Lemma \ref{lem:odd-p-property}}\label{sec:odd-p-property}

We present the argument given in \cite{AC}.
Let $Z_i(G)$ be the upper central series of $G$ and $\gamma_i(G)$ its lower central series. We firstly show that there is an element $t$ of order $p$ such that $t\in Z_2(G)- Z(G)$, i.e., a noncentral element such that for every $g\in G$ the commutator $[g,t]$ is central.

Since $G$ is nonabelian, the set
\[
U=\{j :\, j>1\textrm{, and there exists }x\in Z_{j}(G)- Z_{j-1}(G)\textrm{ of order }p\}
\]
is nonempty. Indeed, otherwise all of the generators of $G$ would be central elements, thus $G$ would be an abelian group. Let $m=\min (U)$.

Now assume that $m>2$. Let $x$ be an element of order $p$ such that $x\in Z_m(G)- Z_{m-1}(G)$. Since $x\notin Z_{m-1}(G)$, there is an element $u\in G$ such that $[x,u]\notin Z_{m-2}(G)$. In particular, since $m>2$, we have that $[x,u]\notin Z(G)$. Let $k$ be the greatest such that there is an element $u\in \gamma_k(G)$ with $[x,u]\notin Z(G)$. 

Let $u$ be such an element. By the definition of lower central series, we have that $[x,u]\in \gamma_{k+1}(G)$. Thus if $[x,[x,u]]\notin Z(G)$, this would contradict the maximality of $k$, so we obtain that $[x,[x,u]]\in Z(G)$. Therefore the group $H=\left\langle x,[x,u]\right\rangle$ has nilpotence class at most two.

Note that $x^u=u^{-1}xu=x[x,u]\in H$. Since in a group of nilpotence class at most two all commutators are central, we have that
$$
\begin{aligned}
[x,u]^p  &= (x^{-1}x^u)^p\\
&= x^{-p} x^{pu} [x^u,x^{-1}]^{\binom{p}{2}}\\
&= [[x,u],x^{-1}]^{\binom{p}{2}} \\
&= [[x,u],x^{-\binom{p}{2}}] =1
\end{aligned}
$$
where we used the following commutator identities: $(gh)^n = g^nh^n[h,g]^{\binom{p}{2}}$ for two elements whose commutator is central, $[gh,k]=[g,k]^h[h,k]$ and $[g^{-1},h]=[h,g]^{g^{-1}}$. In the last line we used the fact that  $\binom{p}{2}$ is divisible by $p$ since $p>2$.
Thus $[x,u]$ is an element of order $p$ in some $Z_j(G)$ with $m>j>1$, contradicting minimality of $m$. So we obtain that $m=2$.  

Now let $t$ be an element of order $p$ such that $Z_2(G)- Z(G)$. Since $t\notin Z(G)$ and $G$ is generated by elements of order $p$, there exists $x$ of order $p$ such that $[x,t]\neq 1$. But in the group $K=\left\langle x,t\right\rangle$ of nilpotence class two we now have
\[
(xt)^p=x^pt^p[t,x]^{\binom{p}{2}}=[t,x^{\binom{p}{2}}]=1.
\]
Therefore $x^{-1}$ and $t$ is the pair of elements satisfying the conclusion.

\end{document}